\documentclass{amsart}
\usepackage{amssymb,amsfonts}
\usepackage{enumerate}
\usepackage{mathrsfs}
\usepackage{graphicx}
\usepackage{mathtools}
\usepackage{leftidx}
\usepackage{tikz-cd}
\usepackage{bm}
\usepackage{url}
\newtheorem{thm}{Theorem}[section]
\newtheorem{cor}[thm]{Corollary}
\newtheorem{prop}[thm]{Proposition}
\newtheorem{lem}[thm]{Lemma}

\theoremstyle{definition}
\newtheorem{defn}[thm]{Definition}

\newtheorem{exmp}[thm]{Example}

\theoremstyle{remark}
\newtheorem{rem}[thm]{Remark}

\makeatletter
\let\c@equation\c@thm
\makeatother
\numberwithin{equation}{section}

\def\bthm{\begin{thm}}
\def\ethm{\end{thm}}
\def\blm{\begin{lem}}
\def\elm{\end{lem}}
\def\bdf{\begin{defn}}
\def\edf{\end{defn}}
\def\bpf{\begin{proof}}
\def\epf{\end{proof}}
\def\bpp{\begin{prop}}
\def\epp{\end{prop}}
\def\bcor{\begin{cor}}
\def\ecor{\end{cor}}
\def\brm{\begin{rem}}
\def\erm{\end{rem}}
\def\beg{\begin{exmp}}
\def\eeg{\end{exmp}}
\def\bN{\mathbb{N}}

\def\bQ{\mathbb{Q}}

\def\bZ{\mathbb{Z}}
\def\cA{\mathcal{A}}
\def\cB{\mathcal{B}}
\def\cC{\mathcal{C}}
\def\cD{\mathcal{D}}

\def\cP{\mathcal{P}}

\def\cR{\mathcal{R}}
\def\cS{\mathcal{S}}

\def\scE{\mathscr{E}}

\def\scO{\mathscr{O}}

\newcommand{\raq}{\,\rightarrow \,}

\newcommand{\rintoq}{\,\hookrightarrow\,}

\newcommand{\xraq}[2][]{\, \xrightarrow[#1]{#2} \,}

\newcommand{\ra}{\rightarrow}

\newcommand{\rinto}{\hookrightarrow}

\newcommand{\ronto}{\twoheadrightarrow}

\newcommand{\xra}[2][]{\xrightarrow[#1]{#2}}

\newcommand{\xronto}[2][]{\xrightarrow[#1]{#2}\mathrel{\mkern-14mu}\rightarrow}

\newcommand{\Mod}{{\rm Mod}}

\newcommand{\Ch}{{\rm Ch}}
\newcommand{\Tot}{{\rm Tot}}
\newcommand{\cone}{{\rm cone}}

\newcommand{\RHomcom}{{\bm R} \underline{{\rm Hom}}}

\newcommand{\HH}{HH}
\newcommand{\HC}{HC}
\newcommand{\HN}{HN}

\newcommand{\CC}{CC}
\newcommand{\CN}{CN}

\newcommand{\Ob}{\rm{Ob}}
\newcommand{\op}{{\rm op}}

\newcommand{\id}{{\rm id}}

\newcommand{\Hom}{{\rm Hom}}
\newcommand{\Homcom}{\underline{{\rm Hom}}}

\newcommand{\coker}{{\rm coker}}

\newcommand{\cHom}{\mathscr{H}\text{\kern -3pt {\calligra\large om}}\,}

\newcommand{\Der}{{\rm Der}}
\newcommand{\Dercom}{\underline{{\rm Der}}}

\newcommand{\Dperf}{\cD_{{\rm perf}}}

\newcommand{\Chdg}{\underline{\Ch}}
\newcommand{\Moddg}{\underline{\Mod}}

\newcommand{\cAe}{\cA^{e}}
\newcommand{\cBe}{\cB^{e}}

\newcommand{\Pie}{\Pi^{e}}

\newcommand{\dgcat}{{\rm dgcat}}

\newcommand{\dgcatflat}{{\rm dgcat}^{{\rm flat}}}

\newcommand{\cof}{{\rm cof}}
\newcommand{\Ho}{{\rm Ho}}

\newcommand{\cRA}{\cR(\cA)}
\newcommand{\cSA}{\cS(\cA)}
\newcommand{\cSB}{\cS(\cB)}

\newcommand{\dtot}{d_{{\rm tot}}}
\newcommand{\tot}{{\rm tot}}

\DeclareMathAlphabet{\mathpzc}{OT1}{pzc}{m}{it}

\title{Relative Calabi-Yau completions}
\author{Wai-Kit Yeung}
\address{Kavli IPMU, The University of Tokyo}
\email{wai-kit.yeung@ipmu.jp}
\begin{document}

\begin{abstract}
We generalize Keller's construction \cite{Kel11} of deformed $n$-Calabi-Yau completions to the relative context.
This gives a construction that extends any given dg functor $F : \cA \ra \cB$ between smooth dg categories
to a dg functor $\widetilde{F} : \widetilde{\cA} \ra \widetilde{\cB}$,
together with a family of deformations of $(\widetilde{\cA},\widetilde{\cB},\widetilde{F})$ parametrized by relative negative cyclic homology classes $\widetilde{\eta} \in \HN_{n-2}(\cB, \cA)$.
We prove that these extensions admit relative $n$-Calabi-Yau structures in the sense of \cite{BD19}. In particular, this proof covers the original absolute case of \cite{Kel11}.
\end{abstract}

\maketitle


\tableofcontents

\section{Introduction}

In \cite{Kel11}, Keller introduced a notion of deformed $n$-Calabi-Yau completions, which constructs a dg algebra $\Pi_n(A,\eta)$ out of the data of a smooth dg algebra $A$ and a Hochschild homology class $\eta \in \HH_{n-2}(A)$. 
As special cases, it recovers the Ginzburg dg algebra \cite{Gin}, as well as the (derived) deformed preprojective algebra \cite{CBH98}. 
Deformed Calabi-Yau completions has found applications in, {\it e.g.}, categorification of cluster algebras. It is also a standard class of examples to consider in quiver representation theory, especially as it relates to Donaldson-Thomas theory.
It is a rich subject because, on the one hand, it is a natural and basic construction when interpreted in a broader picture of formal noncommutative symplectic geometry; on the other hand, it is explicit enough to allow formulaic computations. 

In \cite{Kel11}, it is claimed that the deformed $n$-Calabi-Yau completion is always $n$-Calabi-Yau. This turns out to be incorrect in general. We give a counter-example in Remark \ref{CYcom_not_CY_remark} below. Instead, we formulate and prove a correct version of the statement:

\bthm[= Theorem \ref{main_thm_abs}]  \label{main_thm1_intro}
If the deformation parameter $\eta \in \HH_{n-2}(A)$ has a lift to a negative cyclic class $\widetilde{\eta} \in \HN_{n-2}(A)$, then the deformed $n$-Calabi-Yau completion $\Pi_n(A,\eta)$ has an $n$-Calabi-Yau structure (see Definition \ref{CY_def}).

Moreover, if the deformation parameter is exact, meaning that $\widetilde{\eta}$ is in the image of $B : \HC_{n-3}(A) \ra \HN_{n-2}(A)$, then this $n$-Calabi-Yau structure on $\Pi_n(A,\eta)$ is also exact (see Definition \ref{CY_def}).
\ethm

In fact, we will also extend Keller's construction of deformed $n$-Calabi-Yau completions to the relative context (see Definition \ref{relCYcom_def}), and prove an analogue of Theorem \ref{main_thm1_intro} in this relative case (see Theorem \ref{main_thm_rel} below). Recently, this relative Calabi-Yau completion has found applications in the study of quiver varieties in \cite{BCS}, as well as in the categorification of cluster algebras with frozen variables in \cite{Wu1, Wu2}. 

Now we discuss our proof of Theorem \ref{main_thm1_intro}. 
In the earlier version of this paper, which the author posted on arXiv in 2016, Theorem \ref{main_thm1_intro} was only proved for the case when $A$ is finitely cellular (this covers most of the interesting examples). In the present updated version, we prove it in full generality. While the present version is mostly rewritten from scratch (with a lot of unnecessary material removed), the proof is essentially parallel to the proof in the earlier version for the finitely cellular case, except that at several points, certain equalities are replaced by homotopy equivalences, and one has to keep track of these homotopies throughout various constructions.

In both the finitely cellular and the general case, one encounters the same subtle point that makes the proof of Theorem \ref{main_thm1_intro} non-trivial.
We now explain this subtle point and our solution to it. In particular, our proof will exhibit the interplay between the conceptual and formulaic aspects of (deformed) Calabi-Yau completions that we mentioned above.

Denote by $\Pi = \Pi_n(A,\eta)$ the deformed $n$-Calabi-Yau completion. To show that $\Pi$ is $n$-Calabi-Yau, we need to perform two tasks:
\begin{enumerate}
	\item Construct a negative cyclic class $\widetilde{\omega} \in \HN_n(\Pi)$, whose underlying Hochschild class is denoted as $\omega \in \HH_n(\Pi)$.
	\item Show that the induced map $\omega^{\#} : \Pi^![n] \ra \Pi$ in the derived category $\cD(\Pie)$ of $\Pi$-bimodules is an isomorphism.
\end{enumerate}

Let us explain why it is reasonable to expect an isomorphism $\Pi^![n] \cong \Pi$ in $\cD(\Pie)$, but quite subtle to actually prove it.
From the definition of the deformed Calabi-Yau completion $\Pi$, it is not difficult to see that there is an exact triangle (see Proposition \ref{Omega_Pi_SES} and Corollary \ref{TAM_cofib_cor}):
\begin{equation}  \label{Pi_exact_triag_intro}
	\ldots \raq {\bm L}\iota^e_!( A ) \raq \Pi \raq {\bm L}\iota^e_!(A^!)[n] \raq \ldots 
\end{equation}
where ${\bm L}\iota^e_! : \cD(A^e) \ra \cD(\Pie)$ is the extension functor for bimodules with respect to the natural inclusion $\iota : A \ra \Pi$. Taking the derived dual and shifting by $n$, one obtains an exact triangle
\begin{equation}  \label{Pi_exact_triag_dual_intro}
	\ldots \raq {\bm L}\iota^e_!( A ) \raq \Pi^![n] \raq {\bm L}\iota^e_!(A^!)[n] \raq \ldots 
\end{equation}

Comparing \eqref{Pi_exact_triag_intro} and \eqref{Pi_exact_triag_dual_intro}, it seems reasonable to expect that there is a map of exact triangles from \eqref{Pi_exact_triag_dual_intro} and \eqref{Pi_exact_triag_intro} that are the identity on $\iota^e_!( A )$ and $\iota^e_!(A^!)[n]$, which therefore gives an isomorphism $\Pi^![n] \cong \Pi$ by the $5$-lemma. 

However, to prove that the map $\omega^{\#} : \Pi^![n] \ra \Pi$ induced by a given $\omega \in \HH_n(\Pi)$ indeed satisfies this ``triangularity'', the only way we know is to establish it at the strict level, using explicit models of \eqref{Pi_exact_triag_intro} and \eqref{Pi_exact_triag_dual_intro} in the non-derived category of bimodules.
On the other hand, our construction of $\widetilde{\omega}$, and hence $\omega$, involves the Connes-Tsygan operator on the Hochschild complex. Even in the undeformed case where the construction of $\widetilde{\omega}$ is quite simple, we still find it difficult to prove this ``triangularity'' using the usual definition of the Connes-Tsygan operator.

The major technical innovation of the present paper, which resolves this problem, is to make a systematic use of a certain small model of the negative cyclic complex, given in terms of an $\bN$-graded mixed complex $(X^{\bullet}(\Pi),b,B)$, as developed in \cite{Yeu3} (in fact, our proof of Theorem \ref{main_thm1_intro}, together with some general principles of formal noncommutative algebraic geometry, motivated the constructions in \cite{Yeu3}). The results from \cite{Yeu3} that we need are summarized in Section \ref{sec_small_model} below.

Using this model for the negative cyclic homology, the class $\omega \in \HH_n(\Pi)$ is explicitly modelled by an element $\omega \in X^{(2)}(\Pi)$, which naturally gives a map $\omega^{\#} : \cS(\Pi)^{\vee}[n] \ra \cS(\Pi)$ in the non-derived category $\Mod(\Pie)$ of $\Pi$-bimodules, giving a representative of $\omega^{\#} : \Pi^![n] \ra \Pi$ in $\cD(\Pie)$. Using this explicit strict model, one can verify that the above-mentioned ``triangularity'' holds at the strict level.

In a broader perspective, the employment of the noncommutative calculus from \cite{Yeu3} in our proof of Theorem \ref{main_thm1_intro} shows how this calculus brings together the conceptual and formulaic aspect of formal noncommutative algebraic geometry. This point is further illustrated in \cite{Yeu2}, where the system of noncommutative calculus is further developed.

\vspace{0.2cm}

\textbf{Acknowledgement.} The author thanks Yuri Berest and Alimjon Eshmatov for sharing their initial results that led to the research project in \cite{BEY17, BEY}, from which the present paper stemmed. He also thanks Bernhard Keller for helpful discussions.

\section{Bimodules and duals}

Throughout this paper, we fix a commutative ring $k$ with unit. 
Unadorned tensor product will be understood to be over $k$.
We will work with homological gradings, {\it i.e.}, the differential decreases the degree by $1$.
By a differential graded (dg) category over $k$, we mean a category $\cA$ enriched over the category $\Ch(k)$ of chain complexes of $k$-modules.
We denote the category of all small dg categories over $k$ by $\dgcat_k$.
Basic notions about dg categories, dg functors, modules, {\it etc.} are defined, for example, in \cite{Kel06}.
We recall some of these notions in order to set up convention.

The category $\Ch(k)$ of chain complexes over $k$ can be enriched to a dg category $\Chdg(k)$
where the Hom sets between chain complexes $M$ and $N$ are replaced by Hom complexes $\Homcom_k(M,N)$, 
so that $\Ch(k) = Z_0(\Chdg(k))$.

A right (resp. left) module over a dg category $\cA$ is a dg functor $M : \cA^{\op} \ra \Chdg(k)$ (resp. $M : \cA \ra \Chdg(k)$).
Explicitly, a right module $M$ associates to each object $x\in \cA$ a chain complex $M(x)$, 
together with product maps $M(y) \otimes \cA(x,y) \ra M(x)$ which are associative and unital in the obvious sense.
Unless otherwise stated, a module will always mean a right module.
We denote the category of all right modules over $A$ by $\Mod(\cA)$.

Given $M, N \in \Mod(\cA)$, then a \emph{pre-map of degree $i$} from $M$ to $N$ is a map of graded modules $f : M \ra N[-i]$ over the graded $k$-category $\cA$. Let $\Homcom_{\cA}(M,N)_i$ be the $k$-module of all pre-maps of degree $i$ from $M$ to $N$. Then $\Homcom_{\cA}(M,N)_{\bullet}$ forms a chain complex with $d(f) = d_N \circ f - (-1)^{|f|} f \circ d_M$. Clearly, a map from $M$ to $N$ ({\it i.e., } a morphism in $\Mod(\cA)$) is the same as a degree $0$ pre-map that is closed under $d$. Thus, $\Mod(\cA)$ has a canonical dg enrichment 
$\Moddg(\cA)$ where the Hom sets between dg modules $M,N$ are replaced by Hom complexes $\Homcom_{\cA}(M,N)$.
Denote also by $\Moddg_0(\cA)$ the category of modules and pre-maps of degree $0$.

For any dg categories $\cA$ and $\cB$, we define their tensor product $\cA \otimes \cB$ to be the dg category 
with object set $\Ob(\cA) \times \Ob(\cB)$ and Hom complexes 
\begin{equation*}
	\Homcom_{\cA \otimes \cB}((a,b),(a',b')) := \Homcom_{\cA}(a,a') \otimes \Homcom_{\cB}(b,b')
\end{equation*}
with the obvious composition maps.

In particular, for any dg category $\cA$, we define its \emph{enveloping dg category} 
to be the tensor product $\cAe := \cA \otimes \cA^{\op}$.
Since $\cAe$ is isomorphic to its opposite, left and right modules are the same, and both can be identified with the notion of bimodules:
\bdf 
A \emph{bimodule} $M$ over a dg category $\cA$ associates to each pair $(x,y)\in \Ob(\cA) \times \Ob(\cA)$ 
of objects in $\cA$ a chain complex $M(x,y) \in \cC(k)$, together with maps 
\begin{equation}  \label{bimodule_comp}
	\cA(y_1,y_2) \otimes M(x_2,y_1) \otimes \cA(x_1,x_2) \ra M(x_1,y_2)
\end{equation} 
of chain complexes, which are associative and unital in the obvious sense.
\edf

It is clear that $\cA$ is a bimodule over itself.
For a bimodule $M$, we will often write ${}_{y}M_{x} := M(x,y)$. 
In particular, \eqref{bimodule_comp} can be rewritten as ${}_{y_2}\cA_{y_1} \otimes {}_{y_1}M_{x_2} \otimes {}_{x_2}\cA_{x_1} \ra {}_{y_2}M_{x_1}$. This consistency of object placement in the subscripts often serve as a useful sanity test to ensure that operations are well-defined.
Notice also that, unlike for dg algebras ({\it i.e.,} dg categories with one object), A bimodule $M$ is not automatically a left or right module. Instead, for each $x \in \Ob(\cA)$, there is a left module ${}_{\bullet}M_x$, and for each $y \in \Ob(\cA)$, there is a right module ${}_yM_{\bullet}$.

The same discussion applies to $\cAe$, so that $\cAe$ is itself a $\cAe$-bimodule. To avoid confusion with taking the opposite, we will simply write this $\cAe$-bimodule as $\cA \otimes \cA$. 
Thus, for any pair $(y,x)$, there is a right $\cAe$-module ${}_{(y,x)}(\cAe)_{(\bullet,\bullet)}$, and for any pair $(x',y')$, there is a left $\cAe$-module ${}_{(\bullet,\bullet)}(\cAe)_{(x',y')}$. 
When written in terms of $\cA \otimes \cA$, the right $\cAe$-module ${}_{(y,x)}(\cAe)_{(\bullet,\bullet)}$ translates to the \emph{inner bimodule} 
${}_y\cA_R \otimes {}_L\cA_{x}$ where we use the subscript $L,R$ to indicate the position of left and right action; while the left $\cAe$-module ${}_{(\bullet,\bullet)}(\cAe)_{(x',y')}$ translates to the \emph{outer bimodule} ${}_L\cA_{x'} \otimes {}_{y'}\cA_{R}$.

It is clear that, if $M$ is a right module over a small dg category $\cB$, then there is a left module $\Homcom_{\cB}(M,\cB)$. Applying this to $\cB = \cAe$, we have the following

\bdf
Given a bimodule $M$ over a small dg category $\cA$, its \emph{bimodule dual} is the bimodule $M^{\vee}$ given by
\begin{equation*}
	{}_{\bullet}(M^{\vee})_{\bullet} \, := \, \Homcom_{\cAe}( \, {}_L M_{R} \, , \, {}_{\bullet}\cA_R \otimes {}_L \cA_{\bullet} \, )
\end{equation*}
where we take the Hom complex of $\cA$-bilinear pre-maps from $M$ to $\cA \otimes \cA$, endowed with the inner bimodule structure, which then inherits a bimodule structure from the outer bimodule structure on $\cA \otimes \cA$.
\edf

Given bimodules $M, N$, the complex $M \otimes_{\cAe} N$ is given by
\begin{equation*}
	M \otimes_{\cAe} N \, := \, \bigoplus_{x,y \in \Ob(\cA)} \, {}_y M_x \otimes {}_x N_y \, / \, ( \, (f \xi g) \otimes \eta - (-1)^{|f|(|\xi| + |g| + |\eta|)} \xi \otimes (g \eta f) \,)
\end{equation*}
where we quotient by the $k$-linear span of the displayed relations, for $f \in {}_{x_0}\cA_{x_3}$, $\xi \in {}_{x_3}M_{x_2}$, $g \in {}_{x_2}\cA_{x_1}$, $\eta \in {}_{x_1}N_{x_0}$.

For $N = \cA$, this tensor product will be denoted as $M_{\natural} := M \otimes_{\cAe} \cA$. It can be described as
\begin{equation*}
	M_{\natural}  \, = \, \Bigl( \, \bigoplus_{x \in \Ob(\cA)} M(x,x) \Bigr) / (\xi f - (-1)^{|f||\xi|} f \xi)_{f \in \cA(x,y), \xi \in M(y,x)}
\end{equation*}

Any homogeneous element $\theta \in M \otimes_{\cAe} N$ determines two pre-maps
\begin{equation}  \label{tensor_induce_map}
	\begin{split}
		\theta^{\#}_L \, &: \, M^{\vee} \raq N  \\
		\theta^{\#}_R \, &: \, N^{\vee} \raq M 	
	\end{split}
\end{equation}
defined in the obvious way. These two maps are related by the commutativity of the following diagrams:
\begin{equation}  \label{tensor_induce_map_rel_1}
	\begin{tikzcd}
		M^{\vee} \ar[r, "\theta_L^{\#}"'] \ar[rr, bend left, "(\theta_R^{\#})^{\vee}"] & N \ar[r, "u"'] & N^{\vee\vee} & & N^{\vee} \ar[r, "\theta_R^{\#}"'] \ar[rr, bend left, "(\theta_L^{\#})^{\vee}"] & M \ar[r, "u"'] & M^{\vee\vee}
	\end{tikzcd}
\end{equation}
where $u$ are the canonical maps to double duals.

Denote by $\tau : M \otimes_{\cAe} N \xra{\cong} N \otimes_{\cAe} M$ the switching map. Then we have 
\begin{equation}  \label{tensor_induce_map_rel_2}
(\tau(\theta))^{\#}_R \, = \, \theta^{\#}_L \qquad \text{and} \qquad 
(\tau(\theta))^{\#}_L \, = \, \theta^{\#}_R
\end{equation}
This applies in particular to the case $M = N$, where $\tau$ is the generator of a $\bZ/2$-action.

\bdf
Let $\cA$ be a small dg category and $M$ a right module.
For $x \in \Ob(\cA)$, denote by $h_{x}$ the representable module $h_x = {}_{x}\cA_{\bullet}$. 

The module $M$ is said to be \emph{semi-free} if there is an indexed set $\{x_{\alpha}\}_{\alpha \in S}$ of objects $x_{\alpha} \in \Ob(\cA)$, together with an isomorphism $M \cong \bigoplus_{\alpha \in S} \, h_{x_{\alpha}}$ of the underlying graded modules over the graded category $\cA$ after forgetting the differentials. 
Under the isomorphism $M \cong \bigoplus_{\alpha \in S} \, h_{x_{\alpha}}$, 
the identity elements $1 \in h_{x_{\alpha}}(x_{\alpha})$ is sent to certain elements $\xi_{\alpha} \in M(x_{\alpha})$. The indexed set of elements $\{\xi_{\alpha}\}_{\alpha \in S}$ is called a \emph{basis set} of $M$.

The module $M$ is said to be \emph{cellular} if it is semi-free over a basis set $\{\xi_{\alpha}\}_{\alpha \in S}$, and if $S$ has a filtration $S_1 \subset S_2 \subset \ldots$, such that $S = \bigcup S_i$ and that, 
for each $\beta \in S_i, \, i=1,2,\ldots$, 
the differential $d(\xi_{\beta})$ lies in the graded submodule $\bigoplus_{\alpha \in S_{i-1}} \, h_{x_{\alpha}}$.
Thus, in particular, $d(\xi_{\beta}) = 0$ for all $\beta \in S_1$.

The module $M$ is said to be \emph{finitely cellular} if it is cellular over a finite basis set.
\edf

\bdf
A \emph{graded quiver} consists of a set $\scO$ and a graded set $Q$ ({\it i.e.,} a set $Q$ together with a map of sets $|-| : Q \ra \bZ$), together with two maps $s,t : Q \ra \scO$. We will often simply say that $Q$ is a graded quiver over $\scO$. 

Given a graded quiver $Q$ over $\scO$, then we denote by $T_{\scO}Q$ the graded $k$-category whose object set is $\scO$, and whose morphism space from $x$ to $y$ is the graded $k$-module with a basis consisting of a string of composable arrows (including the empty string).
Clearly, a map ({\it i.e.,} a graded $k$-linear functor) from $T_{\scO}Q$ to another graded $k$-category $\cA$ is the same as a map from $Q$ to the underlying graded quiver of $\cA$.

A dg category $\cA$ is said to be \emph{semi-free} if its underlying graded $k$-category (forgetting the differential) is isomorphic to $T_{\scO}Q$ for some graded quiver $(\scO,Q)$.

A semi-free dg category $\cA \cong T_{\scO}Q$ is said to be \emph{cellular} if the set $Q$ has a filtration $Q_1 \subset Q_2 \subset \ldots$ such that $Q = \bigcup Q_i$ and that, for each $f \in Q_i$, $i=1,2,\ldots$, the corresponding morphism $f \in \cA(s(f),t(f))$ has differential $df$ in the graded subcategory $T_{\scO}(Q_{i-1}) \subset T_{\scO}(Q)$.
It is said to be \emph{finitely cellular} if both $\scO$ and $Q$ are finite.
\edf

\bdf
A dg functor $F : \cA \ra \cB$ is said to be a \emph{quasi-fully faithful} if $F : \cA(x,y) \ra \cB(F(x),F(y))$ is a quasi-isomorphism for any $x,y \in \Ob(\cA)$.
It is said to be a \emph{quasi-equivalence} if it is quasi-fully faithful and $F: H_0(\cA) \ra H_0(\cB)$ is essentially surjective.
\edf

The category $\dgcat_k$ will be assumed to be endowed with the model structure in \cite{Tab05}. In particular, weak equivalences are quasi-equivalences, and cofibrant dg categories are retracts of cellular dg categories. Similarly, $\Mod(\cA)$ will be endowed with the standard model structure where weak equivalences and fibrations are defined pointwise. Cofibrant objects in $\Mod(\cA)$ are retracts of cellular modules.
The derived category $\cD(\cA)$ is the homotopy category of this model category, which is also equivalent to $H_0$ of the full dg subcategory $\Moddg^{\cof}(\cA) \subset \Moddg(\cA)$ consisting of cofibrant objects.
We will denote by $[-] : \Mod(\cA) \ra \cD(\cA)$ and $[-] : \dgcat_k \ra \Ho(\dgcat_k)$ the localization functors to their homotopy categories.

\bdf  \label{lin_cof_flat_def}
We say that $\cA$ is \emph{linearly cofibrant} if each $\cA(x,y) \in \Ch(k)$ is cofibrant for any $x,y \in \scO$.

Recall that a complex $V \in \Ch(k)$ is said to be \emph{h-flat} if the functor $-\otimes V : \Ch(k) \ra \Ch(k)$ preserves quasi-isomorphisms. We say that $\cA \in \dgcat_k$ is \emph{$k$-flat} if $\cA(x,y)$ is h-flat for any $x,y \in \scO$.
We denote by $\dgcatflat_k \subset \dgcat_k$ the full subcategory consisting of $k$-flat small dg categories.

Both conditions are automatic if $k$ is a field. Notice that, if $\cA$ is linearly cofibrant, then $\cRA$ is a cofibrant resolution of $\cA$ as a bimodule; if $\cA$ is $k$-flat, then $\cRA$ is a flat resolution of $\cA$ as a bimodule.
\edf

\brm \label{k_flat_bimodule_remark}
Notions that involve the derived category $\cD(\cAe)$ of bimodules tend to be well-behaved only under the $k$-flatness assumption. Morally speaking, this is because $k$-flatness ensures that $\cAe = \cA \otimes_k \cA$ represents $\cA \otimes_k^{{\bm L}} \cA$.
For a more precise reason, see Lemma \ref{F_shriek_gamma_functorial} below.\erm

\bdf
A module $M \in \Mod(\cA)$ is said to be \emph{perfect} if there exists a finitely cellular module $P \in \Mod(\cA)$ such that $[M] \in \cD(\cA)$ is a retract of $[P] \in \cD(\cA)$ inside the derived category $\cD(\cA)$.
The full subcategory of perfect modules will be denoted as $\Dperf(\cA) \subset \cD(\cA)$. 
In particular, applying this to $\cAe$, one has a notion of \emph{perfect bimodule}, the category of which is denoted as $\Dperf(\cAe) \subset \cD(\cAe)$. 

A module $M \in \Mod(\cA)$ is said to be \emph{strictly perfect} if it is a retract inside $\Mod(\cA)$ of a finitely cellular module. Similarly, we have a notion of strictly perfect bimodules. 

A small dg category $\cA \in \dgcat_k$ is said to be \emph{homotopically finitely presented} if there exists a finitely cellular $\cB \in \dgcat_k$ such that $[\cA] \in \Ho(\dgcat_k)$ is a retract of $[\cB] \in \Ho(\dgcat_k)$ inside the homotopy category $\Ho(\dgcat_k)$.

A small $k$-flat%
\footnote{If $\cA$ is not necessarily $k$-flat, then we say that $\cA$ is smooth if some (hence any) $k$-flat resolution $\cA' \xra{\sim} \cA$ is smooth. See also Remark \ref{k_flat_bimodule_remark}.} dg category $\cA \in \dgcatflat_k$ is said to be \emph{smooth} if $\cA$ is perfect as a bimodule over itself.
\edf

Given $\cA \in \dgcat_k$, denote by $\scO := \Ob(\cA)$. For $n \geq 0$, define $\cR_n(\cA) \in \Mod(\cAe)$ by
\begin{equation*}
	\cR_n(\cA) \, := \, \cA \otimes_{\scO} \stackrel{(n+2)}{\ldots} \otimes_{\scO} \cA
\end{equation*}
where we write ${}_y(M \otimes_{\scO} N)_x := \bigoplus_{z \in \scO} \, ({}_y M_z )\otimes ({}_z N_x)$. 
then there is a resolution of $\cA$ in the abelian category $\Mod(\cAe)$:
\begin{equation}  \label{cR_n_A_complex}
	\ldots \raq \cR_2(\cA) \raq \cR_1(\cA)\raq \cR_0(\cA)
\end{equation}

The direct sum total complex of \eqref{cR_n_A_complex} will be denoted as $\cRA$, which admits a quasi-isomorphism $\cRA \xra{\sim} \cA$ as a bimodule. The bimodule $\cRA$ is called the \emph{bar resolution} of $\cA$.

Given a bimodule $M$, then a \emph{pre-derivation of degree $i$} from $\cA$ to $M$ is a derivation $\cA \ra M[-i]$ to the graded module $M[-i]$ over the graded $k$-category $\cA$. {\it i.e.,} $\delta : \cA \ra M$ is a map of graded $\scO$-bimodules of degree $i$ such that $\delta(f\cdot g) = \delta(f) \cdot g + (-1)^{|f|i} f \delta(g)$ for any $f \in {}_z\cA_{y}$ and $g \in {}_y\cA_{x}$. Let $\Dercom_i(\cA,M)$ be the $k$-module of all pre-derivations of degree $i$, then $\Dercom_{\bullet}(\cA,M)$ forms a chain complex with differential $d(\delta) = d_M \circ \delta - (-1)^{|\delta|} \delta \circ d_{\cA}$. A \emph{derivation} is a pre-derivation of degree $0$ that is closed under this differential: $\Der(\cA,M) := Z_0( \Dercom_{\bullet}(\cA,M) )$.

Denote by $\Omega^1(\cA) \in \Mod(\cAe)$ the bimodule of noncommutative Kahler differentials, characterized by the universal property
\begin{equation*}
	\Dercom(\cA,M) \, \cong \, \Homcom_{\cAe}(\Omega^1(\cA), M)
\end{equation*}
for any $M \in \Mod(\cAe)$. The universal derivation will be denoted as $D : \cA \ra \Omega^1(\cA)$. In particular, $D \in \Der(\cA,\Omega^1(\cA)) = Z_0(\Dercom(\cA, \Omega^1(\cA)))$.

The bimodule $\Omega^1(\cA)$ can be identified as $\Omega^1(\cA) = \coker(\cR_2(\cA) \ra \cR_1(\cA))$, the cokernel of the map in \eqref{cR_n_A_complex}. Hence, there is a short exact sequence
\begin{equation*}
	0 \raq \Omega^1(\cA) \xraq{\alpha} \cA \otimes_{\scO} \cA \xraq{\mu} \cA \raq 0
\end{equation*}
where $\mu$ is the multiplication map, and $\alpha$ is the map of bimodules corresponding to the derivation $\cA \ra \cA \otimes_{\scO} \cA$ given by $f \mapsto f \otimes 1 - 1 \otimes f$.

Notice that $\cA \otimes_{\scO} \cA$ is a free bimodule ${}_{\bullet}\cA \otimes_{\scO} \cA_{\bullet} = \bigoplus_{x \in \scO} \, ({}_{\bullet}\cA_x) \otimes ({}_x\cA_{\bullet})$. It will be convenient to denote the basis element $1_x \otimes 1_x \in {}_x(\cA \otimes_{\scO} \cA)_x$ by $E_x$.
In this notation, the map $\alpha : \Omega^1(\cA) \ra \cA \otimes_{\scO} \cA$ is given by
\begin{equation*}
	\alpha(Df) \, = \, f \cdot E_x - E_y \cdot f  \qquad \text{ for } f \in \cA(x,y)
\end{equation*}

Denote by $\cSA \in \Mod(\cAe)$ the cone
\begin{equation*}
	\cSA \, := \, \cone \, [ \, \Omega^1(\cA) \xraq{\alpha} \cA \otimes_{\scO} \cA \,]
\end{equation*}
called the \emph{short resolution} of $\cA$ (or alternatively the \emph{Cuntz-Quillen resolution}, or the \emph{first order resolution}). Clearly, there is a quasi-isomorphism $\cSA \xronto{\sim} \cA$.

As a graded module, we have $\cSA = (\Omega^1(\cA)[1]) \oplus (\cA \otimes_{\scO} \cA)$. Thus, for any $f \in {}_y\cA_x$, there is an element $sDf \in  {}_y\cSA_x$, where $s$ denotes the homological shift.

If $F : \cA \ra \cB$ is a dg functor, then the restriction and induction functors on (right) modules and bimodules will be denoted by 
\begin{equation}  \label{F_restr_ext}
	\begin{split}
		F_! \, &: \,  \Mod(\cA) \raq \Mod(\cB) \, , \qquad \quad M \mapsto M \otimes_{\cA} \cB \\
		F^* \, &: \, \Mod(\cB) \raq \Mod(\cA) \, , \qquad \quad N \mapsto N|_{\cA} \\
		F^e_! \, &: \,  \Mod(\cAe) \raq \Mod(\cBe) \, , \qquad \quad M \mapsto M \otimes_{\cAe} \cBe \\
		F^{e*} \, &: \, \Mod(\cBe) \raq \Mod(\cAe) \, , \qquad \quad N \mapsto N|_{\cAe} \\
	\end{split}
\end{equation}
More formally, $F^*$ sends a dg functor $N : \cB^{\op} \ra \Chdg(k)$ to the composition $N \circ F$; while $F_!$ is a $\Ch(k)$-enriched left Kan extension. In particular, $F_!$ sends a representable module $h_x$ to the representable module $h_{F(x)}$. The same are true for the bimodule versions.

The functors \eqref{F_restr_ext} form Quillen adjunctions $F_! \dashv F^*$ and $F^e_! \dashv F^{e*}$, whose total derived functors will be denoted as ${\bm L} F_! \dashv F^*$ and ${\bm L} F^e_! \dashv F^{e*}$.

\bdf
A dg functor $F : \cA \ra \cB$ between small dg categories is said to be a \emph{derived equivalence} if the induced functor ${\bm L}F_! : \cD(\cA) \ra \cD(\cB)$ is an equivalence.
\edf

\bpp  \label{derived_equiv_criteria}
Given a dg functor $F : \cA \ra \cB$ between small dg categories.
\begin{enumerate}
	\item The adjunction unit $\id \Rightarrow F^* \circ {\bm L}F_!$ is a natural isomorphism if and only if $F$ is quasi-fully faithful. 
	\item The adjunction counit ${\bm L}F_! \circ F^* \Rightarrow \id$ is a natural isomorphism if and only if for all $b,b' \in \Ob(\cB)$, the canonical map ${}_b\cB \otimes^{{\bm L}}_{\cA} \cB_{b'} \ra {}_b\cB_{b'}$ is a quasi-isomorphism%
	\footnote{If $\cA$ and $\cB$ are $k$-flat, then this map can be regarded as a map $\cB \otimes^{{\bm L}}_{\cA} \cB \ra \cB$ in $\cD(\cBe)$, which may be represented explicitly by the map $\cB \otimes_{\cA} \cRA \otimes_{\cA} \cB \ra \cB$. Namely, the $k$-flatness of $\cA$ and $\cB$ implies that ${}_b\cB \otimes_{\cA} \cRA$ is a flat resolution of $F^*({}_b\cB) \in \Mod(\cA)$, for which we can use to compute ${\bm L}F_! (F^*({}_b\cB))$.}.
	This condition is sometimes called a homological surjection in the literature.
	\item If $F$ is a quasi-equivalence, then it is a derived equivalence.
	\item If $\cA$ and $\cB$ are $k$-flat, and if $F$ is a derived equivalence, then so is $F^e : \cAe \ra \cBe$. Moreover, the canonical map $\gamma_F : {\bm L}F^e_!(\cA) \ra \cB$ is an isomorphism in $\cD(\cBe)$.
\end{enumerate}
\epp

\bpf



(1) and (2) are straightforward. (3) is standard (see, {\it e.g.}, \cite[Proposition 3.9]{KL15}).

For (4), notice that since $\cA$ and $\cB$ are $k$-flat, for any $b_1,b_2,b_3,b_4 \in \Ob(\cB)$, we have
\begin{equation*}
	{}_{(b_1,b_2)}(\cBe) \otimes^{{\bm L}}_{\cAe} (\cBe)_{(b_3,b_4)} \, \simeq \,
	( {}_{b_1} \cB \otimes^{{\bm L}}_{\cA} \cB_{b_3} ) \,  \otimes^{{\bm L}}_{k} \, 	( {}_{b_4} \cB \otimes^{{\bm L}}_{\cA} \cB_{b_2} )
\end{equation*}
Explicitly, by the $k$-flatness of $\cA$ and $\cB$, the left hand side can be represented by ${}_{(b_1,b_2)}(\cBe) \otimes_{\cAe} (\cRA \otimes \cRA) \otimes_{\cAe}  (\cBe)_{(b_3,b_4)}$, which is isomorphic to $( {}_{b_1} \cB \otimes_{\cA} \cRA \otimes_{\cA} \cB_{b_3} ) \,  \otimes_{k} \, 	( {}_{b_4} \cB \otimes_{\cA} \cRA \otimes_{\cA} \cB_{b_2} )$, which represents the right hand side.

The map $\gamma_F$ is by definition the map adjoint to the canonical map $\cA \ra F^{e*}(\cB)$, which is an isomorphism in $\cD(\cAe)$ by (1). Since ${\bm L}F^e_!$ is an equivalence, the corresponding map $\gamma_F$ is also an isomorphism in $\cD(\cBe)$.
Alternatively, one can also show that the map $\gamma_F$ in $\cD(\cBe)$, when evaluated in $(b,b') \in \Ob(\cB)^2$, is precisely the map ${}_b\cB \otimes^{{\bm L}}_{\cA} \cB_{b'} \ra {}_b\cB_{b'}$ in part (2).
\epf

If $\cB$ is a graded $k$-category, and if $f \in \cB(x,y)$ is a homogeneous element, then $\cB[f^{-1}]$ is the graded $k$-category $\cB[f^{-1}] = \cB\langle g \rangle/(gf - 1_x, fg - 1_y)$, where $g$ is an arrow from $y$ to $x$ of degree $-|f|$. More generally, if $S$ is a set of homogeneous morphisms $f_{\alpha} \in \cB(x_{\alpha},y_{\alpha})$, then $\cB[S^{-1}]$ is defined in a similar way.
We say that a dg category $\cA$ is \emph{almost semi-free} if there exists a graded quiver $(\scO,Q)$ and a set $S$ (possibly empty) of homogeneous morphisms in $T_{\scO}(Q)$ such that  $\cA \cong T_{\scO}(Q)[S^{-1}]$ as a graded $k$-category.

We will only need the semi-free case ({\it i.e.,} the case when $S = \emptyset$) of the following Lemma, but we state the more general result because it will be useful for topological examples, since chain dg algebras on based loop spaces have convenient models by almost semi-free dg algebras.

\blm  \label{almost_semifree_Omega}
If $\cA$ is almost semi-free $\cA \cong T_{\scO}(Q)[S^{-1}]$, then the bimodule $\Omega^1(\cA)$ is semi-free over the basis $\{ Df \}_{f \in Q}$. Hence, the bimodule $\cSA$ is semi-free over the basis $\{ sDf \}_{f \in Q} \cup \{E_x\}_{x \in \scO}$.
\elm

\bpf
The statement only depends on the underlying graded $k$-category of $\cA$. For any graded bimodule $M$, we have
\begin{equation}  \label{der_as_map_to_inf_ext}
	\Dercom_0(\cA,M) \, = \, \{ \, \text{graded }k\text{-functor of the form } (\id,\delta) : \cA \ra \cA \oplus M \, \}
\end{equation}
where $\cA \oplus M$ is the infinitesimal extension of the graded $k$-category $\cA$ by the graded bimodule $M$.

Notice that a morphism $(g,\xi) \in (\cA \oplus M)(x,y)$ is invertible in the graded $k$-category $\cA \oplus M$ if and only if $g$ is invertible in the graded $k$-category $\cA$. Denote by $u : T_{\scO}(Q) \ra \cA$ the localization map, then any graded $k$-functor of the form $(u,\widetilde{\delta}) : T_{\scO}(Q) \ra \cA \oplus M$ descends uniquely to the domain $\cA$. Moreover, such a map $\widetilde{\delta}$ is uniquely determined by the map of graded sets $\widetilde{\delta}|_Q$. Hence, the graded bimodule $\Omega^1(\cA)$ that represents the functor $\Dercom_0(\cA,-)$ is free over the basis $\{ Df \}_{f \in Q}$.
\epf

\bcor  \label{cellular_cor_1}
\begin{enumerate}
	\item If $\cA$ is cellular, then $\cSA$ is cellular.
	\item If $\cA$ is finitely cellular, then $\cSA$ is finitely cellular.
	\item If $\cA$ is cofibrant, then $\cSA$ is cofibrant.
\end{enumerate}
\ecor

\bpf
(1) and (2) clearly follows from Lemma \ref{almost_semifree_Omega} (the cellular condition for $\cSA$ follows from the cellular condition for $\cA$). For (3), notice that any dg functor $F : \cA \ra \cB$ gives rise to a canonical map 
$\gamma_F : F^e_!(\cSA) \ra \cSB$.
If $\cA \xra{F} \cB \xra{G} \cA$ is a retract ({\it i.e.,} $G \circ F = \id$), then by functoriality, the composition
\begin{equation}  \label{Omega_retract_functorial}
 \cSA) \, = \, G^e_! F^e_! (\cSA) \xraq{G^e_!(\gamma_F)} G^e_! (\cSB) \xraq{\gamma_G} \cSA
\end{equation}
is the identity, so that $\cSA$ is a retract of $G^e_! (\cSB)$. Thus, if $\cB$ is cellular, then by (1) so is $\cSB$, and hence $G^e_! (\cSB)$, so that $\cSA$ is cofibrant.
\epf

\brm
All three parts of Corollary \ref{cellular_cor_1} are true for $\Omega^1(\cA)$ instead of $\cSA$, with the same proofs.
\erm

\bcor  \label{cellular_cor_2}
If $\cA \in \dgcatflat_k$ is homotopically finitely presented, then it is smooth.
\ecor

\bpf
We repeat the proof of Lemma \ref{cellular_cor_1}(3), except that we replace the functors $F^e_!$ and $G^e_!$ in \eqref{Omega_retract_functorial} by their derived versions ${\bm L}F^e_!$ and ${\bm L}G^e_!$. Lemma \ref{F_shriek_gamma_functorial} below establishes enough functoriality for this argument to work.
\epf

\blm  \label{F_shriek_gamma_functorial}
For any $[\cA] \in \Ho(\dgcat_k)$, denote by $\cD([\cA]^e)$ the derived category $\cD(\cAe)$ for any $k$-flat representative $\cA$ of $[\cA]$. This category has a distinguished object $\cA \in \cD(\cAe)$. The assignment $[\cA] \mapsto (\cD(\cAe),\cA)$ is weakly functorial in the following sense:
for any morphism $F : [\cA] \ra [\cB]$ in $\Ho(\dgcat_k)$, there is a functor
${\bm L} F^e_! : \cD([\cA]^e) \ra \cD([\cB]^e)$, together with $\gamma_F : {\bm L} F^e_!(\cA) \ra \cB$, such that $({\bm L} F^e_!,\gamma_F)$ is functorial in $F$ up to natural isomorphisms. {\it i.e.,} for $[\cA] \xra{F} [\cB] \xra{G} [\cC]$, we have $\phi_{G,F} : {\bm L} (G \circ F)^e_! \xra{\cong} {\bm L} G^e_! \circ {\bm L} F^e_!$, and the composition
\begin{equation*}
	{\bm L} (G \circ F)^e_!(\cA) \xraq[\cong]{\phi_{G,F}} {\bm L} G^e_! {\bm L} F^e_! (\cA) \xraq{ {\bm L} G^e_!(\gamma_F) } {\bm L} G^e_!(\cB) \xraq{\gamma_G} \cC
\end{equation*}
is equal to $\gamma_{G \circ F}$. We call this a weak functoriality of $({\bm L} F^e_!,\gamma_F)$ in $F$ (we do not claim any coherence among the natural transformations $\phi_{G,F}$).
\elm

\bpf
We first show an analogous statement of weak functoriality for $F : \cA \ra \cB$ in $\dgcat_k$. Recall that ${\bm L} F^e_! : \cD(\cAe) \ra \cD(\cBe)$ is left adjoint to $F^{e*} : \cD(\cBe) \ra \cD(\cAe)$, so that we may define $\gamma_F$ as the map adjoint to the canonical map $v_F : \cA \ra F^{e*}(\cB)$. The weak functoriality of $({\bm L} F^e_!,\gamma_F)$ then follows from the weak functoriality of $(F^{e*},v_F)$ by taking left adjoints everywhere.

To obtain the statement for $\Ho(\dgcat_k)$ instead of $\dgcat_k$, we apply Proposition \ref{derived_equiv_criteria}.
More precisely, let $\scE$ be the (big) category whose objects consist of $(\cD,X)$ where $\cD$ is a triangulated category and $X$ an object of $\cD$, 
with Hom-sets given by
\begin{equation*}
	\Hom_{\scE}((\cD,X),(\cD',X')) \, = \, \{ \, (F,u) \, | \,F:\cD \ra \cD' \text{ exact functor, and }  u \in \Hom_{\cD'}(F(X), X') \, \}/_{\sim}
\end{equation*}
where $\sim$ is the equivalence relation that says that $(F,u) \sim (F',u')$ if there exists a natural isomorphism $\phi : F \stackrel{\simeq}{\Rightarrow} F'$ such that $u' \circ \phi_X = u$. The first paragraph of our proof shows that there is a well-defined functor $\dgcat_k \ra \scE$ given by $\cA \mapsto (\cD(\cAe),\cA)$. Proposition \ref{derived_equiv_criteria}(3),(4) implies that this functor sends quasi-equivalences between $k$-flat dg categories to isomorphisms in $\scE$. Since $\Ho(\dgcat_k)$ is equivalent to the localization of $\dgcat_k^{{\rm flat}} \subset \dgcat_k$ at quasi-equivalences, the above functor descends to $\Ho(\dgcat_k) \ra \scE$.
\epf

\brm
It is clear from its proof that this weak functoriality descends to the homotopy category with respect to the model structure \cite{Toe07} on $\dgcat_k$ with derived equivalences as weak equivalences.
\erm

Given a bimodule $M$ over $k$-flat $\cA \in \dgcat_k$, denote by $M^{!}$ its \emph{derived bimodule dual} of a bimodule $M$ by $M^{!}$, {\it i.e.,} it is the object $M^! \in \cD(\cAe)$  represented by $Q^{\vee}$ for a cofibrant replacement $Q \xra{\sim} M$.
%

If $F: \cA \ra \cB$ is a dg functor between small dg categories, then for any $M \in \Mod(\cA)$, there is a canonical map in the derived category of left $\cB$-modules $\cD(\cB^{\op})$:
\begin{equation*}
	\cB \otimes^{{\bm L}}_{\cA} \RHomcom_{\cA}(M,\cA) \raq \RHomcom_{\cA}(M,\cB)
\end{equation*}
which is moreover an isomorphism if $M$ is perfect.

In particular, if $\cA$ and $\cB$ are $k$-flat, then applying this to $F^e : \cAe \ra \cBe$, we have
\blm \label{F_shriek_dual_lemma}
If $F: \cA \ra \cB$ is a dg functor between $k$-flat small dg categories, then for any $M \in \cD(\cAe)$, there is a canonical map in $\cD(\cBe)$:
\begin{equation}  \label{F_shriek_dual_1}
	{\bm L}F^e_!(M^!) \raq {\bm L}F^e_!(M)^!
\end{equation}
which is moreover an isomorphism if $M$ is perfect.
\elm

To be more explicit, suppose that $q : Q \xra{\sim} M$ is a cofibrant replacement, and suppose $p : P \xra{\sim} Q^{\vee}$ is a cofibrant replacement, then the map \eqref{F_shriek_dual_1} is represented by
\begin{equation}  \label{F_shriek_dual_2}
\widetilde{p} \, : \, F^e_!(P) \xraq{F^e_!(p)} F^e_!(Q^{\vee}) \raq F^e_!(Q)^{\vee}
\end{equation}
In particular, if $Q$ is perfect, then $\widetilde{p}$ is a quasi-isomorphism.

Now we discuss duals of semi-free bimodules. Suppose $M \in \Mod(\cAe)$ is semi-free over the set $\{ \xi_i \in M(x_i,y_i) \}_{i \in I}$. Then let $\xi_i^{\vee} \in M^{\vee}(y_i,x_i)$ be the element of degree $-|\xi_i|$ specified by
\begin{equation*}
\xi_i^{\vee}(\xi_j) \, = \, 
\begin{cases*}
1_{x_i} \otimes 1_{y_i} \in {}_{x_i}\cA_{x_i} \otimes {}_{y_i}\cA_{y_i}  & if $i = j$ \\
0 \in {}_{x_i}\cA_{x_j} \otimes {}_{y_j}\cA_{y_i} & if $i \neq j$
\end{cases*}
\end{equation*}

Now we assume that $I$ is finite, say $I = \{1,\ldots,m\}$. Then it is easy to see that $\{\xi_1^{\vee},\ldots,\xi_m^{\vee}\}$ forms a basis set of $M^{\vee}$, so that we have a complete description of $M^{\vee}$ provided that we specify $d(\xi_i^{\vee})$. 
Suppose we have
\begin{equation} \label{basis_diff_xi}
d(\xi_i) = \sum_{j=1}^m f_{ij} \cdot \xi_j \cdot g_{ji} 
\end{equation}
for $f_{ij} \in {}_{y_i}\cA_{y_j}$ and $g_{ji} \in {}_{x_j}\cA_{x_i}$, then the differential on the dual basis is given by 
\begin{equation} \label{dualmod_diff}
d( \xi_j^{\vee} ) = - \sum_{i=1}^m  (-1)^{|f_{ij}|(|\xi_i|+|\xi_j|+|g_{ji}|)} 
g_{ji} \cdot \xi_i^{\vee} \cdot f_{ij}
\end{equation}
Indeed, this follows from $\langle d(\xi_j^{\vee}) ,  \xi_i \rangle = d( \langle \xi_j^{\vee} ,  \xi_i \rangle ) - (-1)^{|\xi_j|} \langle \xi_j^{\vee} ,  d(\xi_i) \rangle $. (Notice that the multiplications by $f_{ij}$ and $g_{ij}$ have switched positions when comparing \eqref{basis_diff_xi} and \eqref{dualmod_diff}. This is because the bimodule structure of $M^{\vee}$ is inherited from the outer multiplication on $\cA \otimes \cA$, while the map $\xi^{\vee} : M \ra {}_{x_i}\cA \otimes \cA_{y_i}$ is $\cA$-bilinear with respect to the inner multiplication on $\cA \otimes \cA$.)

Similarly, one can determine the map $\alpha^{\vee} : N^{\vee} \ra M^{\vee}$ induced by a map 
$\alpha : M \ra N$ in terms of basis elements.
Thus, suppose $M$ has a basis $\{\xi_1,\ldots,\xi_m\}$ and $N$ has a basis $\{ \eta_1,\ldots,\eta_n \}$,
and suppose that the map $\alpha$ is given by 
\[ \alpha(\xi_i) = \sum_{j=1}^n f_{ij} \cdot \eta_j \cdot g_{ji}\]
then the induced map $\alpha^{\vee} : N^{\vee} \ra M^{\vee}$ is given by
\begin{equation}  \label{dualmod_map}
\alpha^{\vee}( \eta_j^{\vee} ) = \sum_{i=1}^m  (-1)^{|f_{ij}|(|\xi_i|+|\xi_j|+|g_{ji}|)} 
g_{ji} \cdot \xi_i^{\vee} \cdot f_{ij} 
\end{equation}

\section{A small model for the negative cyclic complex}  \label{sec_small_model}

Given $\cA \in \dgcat_k$, then for $n \geq 0$, let $X^{(n)}(\cA)$ be the complex defined by
\begin{equation*}
	X^{(n)}(\cA) \, := \, (\cSA \otimes_{\cA} \stackrel{(n)}{\ldots} \otimes_{\cA} \cA)_{\natural}
\end{equation*}
whose differential will be denoted by $b$.

In \cite{Yeu3}, we constructed a structure of $\bN$-graded mixed complex on the complexes $X^{(n)}(\cA)$, in the sense of the following
\bdf
An \emph{$\bN$-graded mixed complex} is a sequence $\{ (C^{(n)},b)\}_{n \geq 0}$ of chain complexes, together with maps $B: C^{(n)} \ra C^{(n+1)}$ of homological degree $+1$, satisfying $B^2 = 0$ and $Bb + bB = 0$.
\edf

Thus, there are maps
\begin{equation}  \label{X_mixed_complex}
	X^{(0)}(\cA) \xraq{B} X^{(1)}(\cA) \xraq{B} X^{(2)}(\cA) \xraq{B} \ldots
\end{equation}
satisfying $B^2 = 0$ and $Bb + bB = 0$.

The map $B : X^{(1)}(\cA) \ra X^{(2)}(\cA)$ will be particularly important to the arguments below, so we give a description.

Recall that, as a graded module, we have $\cSA = (\Omega^1(\cA)[1]) \oplus (\cA \otimes_{\scO} \cA)$. Thus, elements of $X^{(1)}(\cA)$ consists of finite sums of elements of two types: 
\begin{enumerate}
	\item $g \cdot sDf$, for $f \in \cA(x,y)$ and $g \in \cA(y,x)$
	\item $f \cdot E_x$ for $x \in \scO$ and $f \in \cA(x,x)$
\end{enumerate}

Similarly, elements of $X^{(2)}(\cA)$ consists of finite sums of elements of four types:
\begin{enumerate}
	\item $f_4 \cdot sDf_3 \cdot f_2 \cdot sDf_1$, for $f_i \in \cA(x_{i-1}, x_i)$ where $x_0 = x_4$.
	\item $f_3 \cdot sDf_2 \cdot f_1 \cdot E_{x_0}$, for $f_i \in \cA(x_{i-1}, x_i)$ where $x_0 = x_3$.
	\item $E_{x_3} \cdot f_3 \cdot sDf_2 \cdot f_1$, for $f_i \in \cA(x_{i-1}, x_i)$ where $x_0 = x_3$.
	\item $E_x \cdot f \cdot E_y \cdot g$, for $f \in \cA(y,x)$ and $g \in \cA(x,y)$.
\end{enumerate}

The following description of $B : X^{(1)}(\cA) \ra X^{(2)}(\cA)$ follows directly from its definition:
\blm  \label{B_on_X1}
On the elements of type (1) in $X^{(1)}(\cA)$, we have 
\begin{equation*}
	B(g \cdot sDf) \, = \, sDg \cdot sDf + (-1)^{(|f|+1)(|g|+1)} sDf \cdot sDg 
\end{equation*}
which is a sum of elements of type (1) in $X^{(2)}(\cA)$.

On the elements of type (2) in $X^{(1)}(\cA)$, we have 
\begin{equation*}
	B(f \cdot E_x) \, = \, sDf \cdot E_x + E_x \cdot sDf
\end{equation*}
which is a sum of elements of type (2) and (3) in $X^{(2)}(\cA)$.
\elm

The map $B : X^{(0)}(\cA) \ra X^{(1)}(\cA)$ can also be described in these terms. Notice that by definition we have $X^{(0)}(\cA) = \cA_{\natural}$.
\blm  \label{B_on_X0}
The map $B : X^{(0)}(\cA) \ra X^{(1)}(\cA)$ sends $f \in \cA_{\natural}$ to the element $sDf$ of type (1) in $X^{(1)}(\cA)$.
\elm

Denote by $C^H(\cA)$ the Hochschild complex $C^H(\cA) := \cRA_{\natural}$, whose differential will be denoted as $b$. Denote by $B : C^H(\cA) \ra C^H(\cA)$ the Connes-Tsygan operator, so that $(C^H(\cA),b,B)$ is a mixed complex.

\bdf
We say that $\cA \in \dgcat_k$ is \emph{almost cofibrant} if it is linearly cofibrant and $\cSA \in \Mod(\cAe)$ is cofibrant.
By Corollary \ref{cellular_cor_1}, if $\cA$ is cofibrant, then it is almost cofibrant.
\edf

We recall the following result from \cite{Yeu3}:

\bthm  \label{X_complex_thm_1}
Suppose that $\cA$ is almost cofibrant.
Then there is a zig-zag of quasi-isomorphisms relating the $\bN_{\geq 1}$-graded part of $(X^{(n)}(\cA),b,B)$ and $(C^H(\cA),b,B)$. More precisely, there is a commutative diagram
\begin{equation}  \label{X_complex_diag_1}
	\begin{tikzcd}
		C^H(\cA) \ar[r, "B"] \ar[d, "\simeq"] & C^H(\cA) \ar[r, "B"] \ar[d, "\simeq"] & C^H(\cA) \ar[r, "B"] \ar[d, "\simeq"] & \ldots \\
		X(\cA) \ar[r, "B"] & X(\cA) \ar[r, "B"] & X(\cA) \ar[r, "B"] & \ldots \\
		X^{(1)}(\cA) \ar[r, "B"] \ar[u, "\simeq"'] & X^{(2)}(\cA) \ar[r, "B"] \ar[u, "\simeq"'] & X^{(3)}(\cA) \ar[r, "B"] \ar[u, "\simeq"'] & \ldots
	\end{tikzcd}
\end{equation}
where the middle row is another $\bN_{\geq 1}$-graded mixed complex, and all the vertical maps are quasi-isomorphisms of chain complexes (with respect to the $b$-differentials).
\ethm

\bdf
Given an $\bN$-graded mixed complex $(C^{\bullet},b,B)$, its \emph{(direct product) total complex} is the complex
\begin{equation*}
	C^{\tot} \, := \, \prod_{n \geq 0} \, C^{(n)} \cdot u^n \, , \qquad \qquad \dtot = b + uB
\end{equation*}
where $u$ is a variable of degree $-2$. It comes with a filtration
\begin{equation*}
	F^r C^{\tot} \, := \, \prod_{n \geq r} \, C^{(n)} \cdot u^n \, , \qquad \qquad \dtot = b + uB
\end{equation*}
\edf

A standard spectral sequence argument shows that a quasi-isomorphism of $\bN$-graded mixed complexes ({\it i.e.,} each $(C^{(n)},b) \ra (C^{\prime (n)},b')$ is a quasi-isomorphism) induces a quasi-isomorphism on the direct product total complexes. 

\bdf
Given a $k$-flat $\cA \in \dgcat_k$, then its \emph{negative cyclic complex} is the complex
\begin{equation*}
	\CN(\cA) \, := \, \prod_{n \geq 0} \, C^{H}(\cA) \cdot u^n \, , \qquad \qquad \dtot = b + uB
\end{equation*}
whose homology is denoted as $\HN_{\bullet}(\cA) := H_{\bullet}(\CN(\cA))$, and is known as the \emph{negative cyclic homology} of $\cA$.
The map $\CN(\cA) \ra C^H(\cA)$ that projects to the $u^0$-component is a map of complexes, and induces a map $h : \HN(\cA) \ra \HH(\cA)$ at the level of homology.
\edf

Notice that the negative cyclic complex is the direct product total complex of the top row of \eqref{X_complex_diag_1}.
By Theorem \ref{X_complex_thm_1}, it is therefore quasi-isomorphic to the direct product total complex of the bottom row of \eqref{X_complex_diag_1}. In other words, we have

\bthm  \label{X_complex_thm_2}
Suppose that $\cA$ is almost cofibrant.
Then for each $r \geq 1$, there is a zig-zag of quasi-isomorphisms relating 
$F^r X^{\tot}(\cA)$ and $\CN(\cA) \cdot u^r = \CN(\cA)[-2r]$.
\ethm

\brm
The entire complex $X^{\tot}(\cA)$ is ``almost'' quasi-isomorphic to the periodic cyclic complex, in the sense that the reduced version of $X^{\tot}(\cA)$ is quasi-isomorphic to the reduced periodic cyclic complex. See \cite{Yeu3} for details.
\erm


The canonical map $B:\HC_{n-1}(\cA) \ra \HN_n(\cA)$ can also be represented by the $\bN$-graded mixed complex $(X^{\bullet}(\cA),b,B)$, or rather its reduced version. We discuss this relation now. This will be used later in the discussion of deformed Calabi-Yau completions with exact deformation parameter.
Readers not interested in exact Calabi-Yau structures may safely skip the rest of this section.

Let $k\scO$ be the dg category with object set $\scO$ and with $k\scO(x,x) = k$ and $k\scO(x,y) = 0$ if $x \neq y$. There is an obvious dg functor $k\scO \ra \cA$.
This induced a map
$\bigoplus_{x \in \scO} k \ra \cA_{\natural}$. Let
\begin{equation*}
\overline{\cA_{\natural}} \, := \, \cone \, [ \,  \bigoplus_{x \in \scO} k \raq \cA_{\natural} \, ]
\end{equation*}

By the description of $B : \cA_{\natural} = X^{(0)}(\cA) \ra X^{(1)}(\cA)$ in Lemma \ref{B_on_X0}, it clearly vanishes on the image of $\bigoplus_{x \in \scO} k \ra \cA_{\natural}$, so that it factors through a map $\overline{B} : \overline{\cA_{\natural}} \ra X^{(1)}(\cA)$, which sends $\bigoplus_{x \in \scO} k$ to zero. 
We modify the $\bN$-graded mixed complex $(X^{\bullet}(\cA),b,B)$ by changing $X^{(0)}(\cA)$ to $\overline{X^{(0)}}(\cA) := \overline{\cA_{\natural}}$:
\begin{equation*}
X^{\bullet}_{{\rm sr}}(\cA) \, := \, [ \, \overline{X^{(0)}}(\cA) \xraq{\overline{B}} X^{(1)}(\cA) \xraq{B} X^{(2)}(\cA) \xraq{B} \ldots \,]
\end{equation*}
where the subscript ``sr'' refers to it being a ``semi-reduced'' version (where we reduce only $X^{(0)}(\cA)$ and keep $X^{(i)}(\cA)$ unchanged for $i > 0$).

For any $r > 0$, consider 
\begin{equation*}
\Tot \, X^{<r}_{{\rm sr}}(\cA) \,:= \, X^{\tot}_{{\rm sr}}(\cA) / F^r X^{\tot}_{{\rm sr}}(\cA) \, = \, \bigoplus_{n = 0}^{r-1} \, X_{{\rm sr}}^{(n)}(\cA) \cdot u^n \, , \qquad d_{\tot} = b + uB
\end{equation*}
which comes with an obvious map $\overline{uB} : \Tot \, X^{<r}_{{\rm sr}}(\cA) \ra F^r X^{\tot}_{{\rm sr}}(\cA)[1]$, such that $X^{\tot}_{{\rm sr}}(\cA)$ is the cocone of this map $\overline{uB}$.

On the other hand, the dg functor $k\scO \ra \cA$ also induces $\bigoplus_{x \in \scO} \CC(k) \ra \CC(\cA)$. Let
\begin{equation*}
\overline{\CC}(\cA) \, := \, \cone \, [ \,  \bigoplus_{x \in \scO} \CC(k) \raq \CC(\cA)  \, ]
\end{equation*}
Then the map $\overline{B} : \CC(\cA) \ra \CN(\cA)[-1]$ also vanishes on the image of $\bigoplus_{x \in \scO} \CC(k) \ra \CC(\cA)$, so that it factors through a map $\overline{B} : \overline{\CC}(\cA) \ra \CN(\cA)[-1]$ which sends $\bigoplus_{x \in \scO} \CC(k)$ to zero.

The following result combines the Feigin-Tsygan theorem \cite{FT85} with Theorem \ref{X_complex_thm_2} (see \cite{Yeu3} for details):

\bthm  \label{X_complex_thm_3}
Assume that $\cA$ is cofibrant and $\bQ \subset k$, then there is a zig-zag of quasi-isomorphisms between the map 
$\overline{uB} : \Tot \, X^{<r}_{{\rm sr}}(\cA) \ra F^r X^{\tot}_{{\rm sr}}(\cA)[1]$
and the map
$\overline{B} : \overline{\CC}(\cA)[-2r+2] \ra \CN(\cA)[-2r+1]$.
\ethm

\brm
In Sections \ref{sec_main_thm_abs} and \ref{sec_main_thm_rel}, when we discuss the deformed Calabi-Yau completion $\Pi = \Pi_n(\cA,\eta)$, it will be convenient to apply Theorem \ref{X_complex_thm_2} in two ways: we use the case $r=1$ to represent the deformation parameter $\widetilde{\eta} \in \HN_{n-2}(\cA)$, and the case $r=2$ to represent the canonical Calabi-Yau structure $\widetilde{\omega} \in \HN_n(\Pi)$ on $\Pi$.

Similarly, in the exact case, we also apply Theorem \ref{X_complex_thm_3} in two ways: the case $r=1$ for the deformation parameter, and the case $r=2$ for the resulting exact Calabi-Yau structure on $\Pi$.
\erm

\section{Calabi-Yau structures}  \label{sec_CY_str}

From now on, we will assume that $\cA$ is $k$-flat, so that $\cD(\cAe)$ is well-behaved (see Remark \ref{k_flat_bimodule_remark}).
By choosing cofibrant resolutions, we see that the derived analogue of \eqref{tensor_induce_map} holds. 
Namely, any element $\theta \in H_0(M \otimes_{\cAe}^{{\bm L}} N)$ determines two maps in the derived category $\cD(\cAe)$:
\begin{equation}   \label{theta_LR_derived}
	\begin{split}
		\theta^{\#}_L \, &: \, M^{!} \raq N  \\
		\theta^{\#}_R \, &: \, N^{!} \raq M 	
	\end{split}
\end{equation}
which are related to each other by the derived analogues of \eqref{tensor_induce_map_rel_1} and \eqref{tensor_induce_map_rel_2}.

In particular, in the case $M = N = \cA$, we see that any element in $\theta \in \HH_n(\cA)$ determines two maps in the derived category $\cD(\cAe)$:
\begin{equation}  \label{theta_LR_Hoch}
	\begin{split}
		\theta^{\#}_L \, &: \, \cA^{!}[n] \raq \cA  \\
		\theta^{\#}_R \, &: \, \cA^{!}[n] \raq \cA 	
	\end{split}
\end{equation}

We recall the following result from \cite[Proposition 14.1]{VdB15} (see also \cite{Yeu3}):

\bpp
The two maps in \eqref{theta_LR_Hoch} are equal: $\theta^{\#}_L = \theta^{\#}_R$. We will denote either by $\theta^{\#}$.
\epp

\bdf  \label{CY_def}
Given a smooth $k$-flat $\cA \in \dgcat_k$, then an \emph{$n$-Calabi-Yau structure} is a negative cyclic class $\widetilde{\omega} \in \HN_n(\cA)$ whose underlying Hochschild class $\omega := h(\widetilde{\omega}) \in \HH_n(\cA)$ is \emph{non-degenerate} in the sense that the induced map $\omega^{\#} : \cA^![n] \xra{\cong} \cA$ is an isomorphism in $\cD(\cAe)$.

We say that an $n$-Calabi-Yau structure is \emph{exact} if the negative cyclic class $\widetilde{\omega} \in \HN_n(\cA)$ is in the image of the canonical map $B : \HC_{n-1}(\cA) \ra \HN_n(\cA)$.
\edf

Now we discuss the notion of relative Calabi-Yau structures. Given a dg functor $F : \cA \ra \cB$ between small $k$-flat dg categories, the \emph{relative Hochschild complex} and the \emph{relative negative cyclic complex} are the complexes
\begin{equation*}
\begin{split}
	C^H(\cB,\cA) \, &:= \, \cone \, [ \, C^H(\cA) \raq C^H(\cB) \, ] \\
	\CN(\cB,\cA) \, &:= \, \cone \, [ \, \CN(\cA) \raq \CN(\cB) \, ]
\end{split}
\end{equation*} 
whose homologies are denoted as $\HH_{\bullet}(\cB,\cA)$ and $\HN_{\bullet}(\cB,\cA)$ respectively. Clearly, in this relative context, there is still a map $h : \HN_{\bullet}(\cB,\cA) \ra \HH_{\bullet}(\cB,\cA)$.

Notice that there are three obvious maps
\begin{equation}  \label{rel_Hoch_induce_map_tensor}
	\begin{split}
		\cone [\cA \otimes_{\cAe}^{{\bm L}} \cA \ra \cB \otimes_{\cBe}^{{\bm L}} \cB ] & \raq \cone [ \cB \otimes_{\cAe}^{{\bm L}} \cA \ra \cB \otimes_{\cBe}^{{\bm L}} \cB ] \, = \, \cB \otimes_{\cBe}^{{\bm L}} \cone [ {\bm L}F^e_!(\cA) \ra \cB ] \\
		\cone [\cA \otimes_{\cAe}^{{\bm L}} \cA \ra \cB \otimes_{\cBe}^{{\bm L}} \cB ] & \raq \cA \otimes_{\cAe}^{{\bm L}} \cA[1] \raq {\bm L}F^e_!(\cA) \otimes_{\cBe}^{{\bm L}} {\bm L}F^e_!(\cA)[1]  \\
		\cone [\cA \otimes_{\cAe}^{{\bm L}} \cA \ra \cB \otimes_{\cBe}^{{\bm L}} \cB ] & \raq \cone [ \cA \otimes_{\cAe}^{{\bm L}} \cB \ra \cB \otimes_{\cBe}^{{\bm L}} \cB ] \, = \, \cone [ {\bm L}F^e_!(\cA) \ra \cB ] \otimes_{\cBe}^{{\bm L}} \cB
	\end{split}
\end{equation}


Thus, any element $\omega \in \HH_{n}(\cB,\cA)$ determines three elements:
\begin{enumerate}
	\item $\omega' \in H_n( \, \cB \otimes_{\cBe}^{{\bm L}} \cone \, [ \, {\bm L}F^e_!(\cA) \raq \cB \,] \, )$
	\item $\omega'' \in H_{n-1}( \, {\bm L}F^e_!(\cA) \otimes_{\cBe}^{{\bm L}} {\bm L}F^e_!(\cA) \, )$
	\item $\omega''' \in H_n( \, \cone \, [ \, {\bm L}F^e_!(\cA) \raq \cB \,] \otimes_{\cBe}^{{\bm L}} \cB \,)$
\end{enumerate}
which in turn determines the three vertical maps in the following diagram (see \eqref{theta_LR_derived}):
\begin{equation}  \label{rel_Hoch_induce_map_1}
	\begin{tikzcd}
		\cB^![n-1] \ar[r] \ar[d, "(\omega')^{\#}_L"]  & {\bm L}F^e_!(\cA)^![n-1] \ar[r]  \ar[d, "(\omega'')^{\#}_L"] &  \cone [ \cB^! \ra {\bm L}F^e_!(\cA)^! ][n-1]  \ar[d, "(\omega''')^{\#}_L"]\\
		\cone [ {\bm L}F^e_!(\cA) \ra \cB ][-1] \ar[r] & {\bm L}F^e_!(\cA) \ar[r] & \cB 
	\end{tikzcd}
\end{equation}

Notice that each row is part of an exact triangle in $\cD(\cBe)$. In fact, the vertical maps give a map of exact triangle:

\bpp  \label{rel_Hoch_induce_map_of_triag_prop}
The vertical maps in \eqref{rel_Hoch_induce_map_1} is a map of exact triangles in $\cD(\cBe)$.
\epp

\bpf
For this verification, it will be convenient to choose a strict model. Choose quasi-isomorphisms $M_{\cA} \xra{\sim} \cA$ and $N_{\cA} \xra{\sim} \cA$ in $\Mod(\cAe)$, and $M_{\cB} \xra{\sim} \cB$ and $N_{\cB} \xra{\sim} \cB$ in $\Mod(\cBe)$, where $M_{\cA}$, $N_{\cA}$, $M_{\cB}$, $N_{\cB}$ are cofibrant.
Choose maps $\gamma_M : F^e_!(M_{\cA}) \ra M_{\cB}$ and $\gamma_N : F^e_!(N_{\cA}) \ra N_{\cB}$ that represents $\gamma_F: {\bm L}F^e_!(\cA) \ra \cB$.
(Of course, we could choose $M_{\cA} = N_{\cA}$, $M_{\cB} = N_{\cB}$ and $\gamma_M = \gamma_N$, but we will write them separately so that it will be clear below which map we are referring to.)

The maps \eqref{rel_Hoch_induce_map_tensor} then admit a strict model:
%
\begin{equation*}  
	\begin{split}
		\cone [M_{\cA} \otimes_{\cAe} N_{\cA} \ra M_{\cB} \otimes_{\cBe} N_{\cB} ] & \raq \cone [ M_{\cB} \otimes_{\cAe} N_{\cA} \ra M_{\cB} \otimes_{\cBe} N_{\cB} ] \, = \, M_{\cB} \otimes_{\cBe} \cone [ F^e_!(N_{\cA}) \ra N_{\cB} ] \\
		\cone [M_{\cA} \otimes_{\cAe} N_{\cA} \ra M_{\cB} \otimes_{\cBe} N_{\cB} ] & \raq M_{\cA} \otimes_{\cAe} N_{\cA}[1] \raq F^e_!(M_{\cA}) \otimes_{\cBe} F^e_!(N_{\cA})[1]  \\
		\cone [M_{\cA} \otimes_{\cAe} N_{\cA} \ra M_{\cB} \otimes_{\cBe} N_{\cB} ] & \raq \cone [ M_{\cA} \otimes_{\cAe} N_{\cB} \ra M_{\cB} \otimes_{\cBe} N_{\cB} ] \, = \, \cone [ F^e_!(M_{\cA}) \ra M_{\cB} ] \otimes_{\cBe} N_{\cB}
	\end{split}
\end{equation*}

Thus, any closed element $\omega \in Z_n( \cone[ M_{\cA} \otimes_{\cAe} N_{\cA} \ra M_{\cB} \otimes_{\cBe} N_{\cB} ] )$ determines the following vertical maps in $\Mod(\cBe)$:
\begin{equation}  \label{rel_Hoch_induce_map_2}
	\begin{tikzcd}
		M_{\cB}^{\vee}[n-1] \ar[r] \ar[d, "(\omega')^{\#}_L"]  & F^e_!(M_{\cA})^{\vee}[n-1] \ar[r]  \ar[d, "(\omega'')^{\#}_L"] &  \cone [ M_{\cB}^{\vee} \ra F^e_!(M_{\cA})^{\vee} ][n-1]  \ar[d, "(\omega''')^{\#}_L"]\\
		\cone [ F^e_!(N_{\cA}) \ra N_{\cB} ][-1] \ar[r] & F^e_!(N_{\cA}) \ar[r] & N_{\cB} 
	\end{tikzcd}
\end{equation}

It is easy to see that both squares commute (strictly). To complete the proof, we need to show that the next square also commute (up to homotopy). This follows from Lemma \ref{cone_three_squares_lemma}.
\epf

\blm  \label{cone_three_squares_lemma}
Denote by $\iota : X \ra \cone[X \ra Y]$ and $\pi : \cone[X\ra Y] \ra Y[1]$ the canonical maps.
Given a diagram
\begin{equation}  \label{cone_three_squares}
	\begin{tikzcd}
		V_1 \ar[r, "\alpha"] \ar[d, "(f_{11}{,}s^{-1}f_{21})"'] \ar[rd, phantom, "(1)" description] & V_2  \ar[r, "\iota"] \ar[d, "f_{12}"] \ar[rd, phantom, "(2)" description] & \cone[V_1 \xra{\alpha} V_2] \ar[r, "-\pi"] \ar[d, "(f_{21}'s^{-1}{,} f_{22})"] & V_1[1] \ar[d, "(f_{11}{[1]}{,}f_{21}s^{-1})"] \ar[d, shift right = 15, phantom, "(3)"] \\
		\cone[W_1 \xra{\beta} W_2][-1]  \ar[r, "\pi{[-1]}"]  & W_1 \ar[r, "\beta"] & W_2 \ar[r, "\iota"] & \cone[W_1 \xra{\beta} W_2]
	\end{tikzcd}
\end{equation}

Suppose that the square (1) and (2) commutes ({\it i.e.,} $f_{12} \circ \alpha = f_{11}$ and $f_{22} = \beta \circ f_{12}$), and suppose that $f_{21}' = -f_{21}$, then the square (3) commutes up to homotopy.
\elm

\bpf
Since the square (2) commutes, there is an induced map on the cone of the horizontal maps in this square, which is given by
\begin{equation}  \label{cone_of_second_square}
(\, f'_{21}s^{-1} \, , \, f_{12}[1] \, , \, f_{22} \,) \, : \, \cone[V_1 \oplus V_2 \xra{(\alpha , \id)} V_2] \raq \cone[W_1 \xra{\beta} W_2]
\end{equation}
Thus, it suffices to notice that the composition
\begin{equation*}
	V_1[1] \xraq{(-\id, \alpha[1], 0)} \cone[V_1 \oplus V_2 \xra{(\alpha , \id)} V_2] \xraq{\eqref{cone_of_second_square}} \cone[W_1 \xra{\beta} W_2]
\end{equation*}
is precisely the rightmost vertical map $(f_{11}{[1]}{,}f_{21}s^{-1})$ in \eqref{cone_three_squares}.
\epf

\brm  \label{non_cof_M_N_remark}
In the proof of Proposition \ref{rel_Hoch_induce_map_of_triag_prop}, we required that $M_{\cA}$, $M_{\cB}$, $N_{\cA}$ and $N_{\cB}$ are cofibrant replacements of $\cA$ and $\cB$, in order to guarantee that \eqref{rel_Hoch_induce_map_2} represents \eqref{rel_Hoch_induce_map_1}. However, our proof actually shows that, given any ($M_{\cA}$, $M_{\cB}$, $N_{\cA}$, $N_{\cB}$, $\gamma_M$, $\gamma_N$), then any closed element $\omega \in Z_n( \cone[ M_{\cA} \otimes_{\cAe} N_{\cA} \ra M_{\cB} \otimes_{\cBe} N_{\cB} ] )$ determines a commutative diagram \eqref{rel_Hoch_induce_map_2}, such that the connecting square is also commutative up to homotopy. 
\erm

\bdf  \label{rel_CY_def}
Given a dg functor $F : \cA \ra \cB$ between smooth $k$-flat small dg categories. Then a \emph{relative $n$-Calabi-Yau structure} on $F$ is a relative negative cyclic class $\widetilde{\omega} \in \HN_n(\cB,\cA)$ whose underlying relative Hochschild class $\omega \in \HH_n(\cB,\cA)$ is \emph{non-degenerate} in the sense that
\begin{enumerate}
	\item Denote by $\omega_{\cA} \in \HH_{n-1}(\cA)$ the image of $\omega$ under the connecting homomorphism $\HH_n(\cB,\cA) \ra \HH_{n-1}(\cA)$, then $\omega_{\cA}$ is non-degenerate in the sense of Definition \ref{CY_def}.
	\item The three vertical maps in \eqref{rel_Hoch_induce_map_1} are isomorphisms in $\cD(\cBe)$.
\end{enumerate}

We say that a relative $n$-Calabi-Yau structure is \emph{exact} if the relative negative cyclic class $\widetilde{\omega} \in \HN_n(\cB,\cA)$ is in the image of the canonical map $B : \HC_{n-1}(\cB,\cA) \ra \HN_n(\cB,\cA)$.
\edf

\brm
In \cite{BD19}, the notion of non-degeneracy does not require condition (1). We add this condition because we think that it is a natural requirement.
\erm

There are some redundancies in the non-degeneracy condition. First, notice that since $\cA$ is assumed to be smooth, by Lemma \ref{F_shriek_dual_lemma} there is a canonical isomorphism ${\bm L}F^e_!(\cA^!) \xra{\cong} {\bm L}F^e_!(\cA)^!$. It can be easily checked that, under this isomorphism, the map $(\omega'')^{\#}_L : {\bm L}F^e_!(\cA)^! \ra {\bm L}F^e_!(\cA)$ in the middle of \eqref{rel_Hoch_induce_map_1} is given by ${\bm L}F^e_!((\omega_{\cA})^{\#}_L)$, where $(\omega_{\cA})^{\#}_L : \cA^![n-1] \ra \cA^!$ is the map induced by $\omega_{\cA} \in \HH_{n-1}(\cA)$.
Thus, the assumption (1) that $\omega_{\cA}$ is non-degenerate already implies that the middle vertical map of \eqref{rel_Hoch_induce_map_1} is an isomorphism.
Then, by Proposition \ref{rel_Hoch_induce_map_of_triag_prop}, it suffices to verify that either of the other vertical maps is an isomorphism.

\section{Calabi-Yau completions}

We will continue to assume that $\cA \in \dgcat_k$ is $k$-flat.
Besides ensuring that $\cD(\cAe)$ is well-behaved (see Remark \ref{k_flat_bimodule_remark}), this also ensures that any cofibrant bimodule $P \in \Mod(\cAe)$ is $\cA$-flat on the left, {\it i.e.,} the functor $-\otimes_{\cA} P : \Mod(\cA) \ra \Mod(\cA)$ preserves quasi-isomorphisms, so that the functor $-\otimes_{\cA}^{{\bm L}} M$ can be represented by any cofibrant resolution $P \xra{\sim} M$ in $\Mod(\cAe)$. This will be important for the construction that follows.

For $M \in \Mod(\cAe)$, denote by $T_{\cA}M$ the dg category with the same objects as $\cA$, and with Hom complexes given by
\begin{equation*}
	T_{\cA}(M) \, := \, (\cA) \, \oplus \,  (M)  \, \, \oplus \, (M \otimes_{\cA} M) \, \oplus \, (M \otimes_{\cA} M \otimes_{\cA} M)  \, \oplus \, \ldots 
\end{equation*}
with the obvious composition maps.

Given a map $\eta : M[-1] \ra \cA$ of bimodules, let $\delta_{\eta} \in \Dercom_{-1}(T_{\cA}(M),T_{\cA}(M))$ be the unique derivation of degree $-1$ such that $(\delta_{\eta})|_{\cA} = 0$ and $(\delta_{\eta})|_M = \eta : M \ra \cA$. Denote by $d_0$ the intrinsic differential of $T_{\cA}(M)$, then since $\eta$ is a map ({\it i.e.,} not just a pre-map), we have $d_0 \delta_{\eta} + \delta_{\eta} d_0 = 0$. By checking on the generators, we also have $\delta_{\eta}^2 = 0$. Thus, if we let $d_{\eta} = d_0 + \delta_{\eta}$, then we have $d_{\eta}^2 = 0$, so that  $(T_{\cA}(M),d_{\eta})$ is a dg category.

Denote by $\Pi$ the dg category $(T_{\cA}(M),d_{\eta})$, and by $\iota : \cA \ra \Pi$ the canonical inclusion.

\bpp  \label{Omega_Pi_SES}
There are canonical exact sequences in the abelian category $\Mod(\Pie)$:
\begin{equation}  \label{TAM_Omega_SES}
	\iota^e_! \, \Omega^1(\cA) \raq \Omega^1(\Pi) \raq \iota^e_! \, M \raq 0
\end{equation}
\begin{equation}   \label{TAM_S_SES}
	\iota^e_! \, \cS(\cA) \raq \cS(\Pi) \raq \iota^e_! \, M[1] \raq 0
\end{equation}
Moreover, if $M \in \Mod(\cAe)$ is cofibrant, then both sequences are left exact ({\it i.e.,} the leftmost arrow is injective).
\epp

\bpf
It suffices to prove the statements for \eqref{TAM_Omega_SES}, since the ones for \eqref{TAM_S_SES} follow. Given any $N \in \Mod(\Pie)$, then by using an argument similar to \eqref{der_as_map_to_inf_ext} (with $(\cA,M)$ in \eqref{der_as_map_to_inf_ext} replaced by $(\Pi,N)$), one can show that giving a graded derivation $\delta : \Pi \ra N$ is the same as giving a graded derivation $\delta_{\cA} :\cA \ra N$ together with a map of $\scO$-bimodules $\delta_M : M \ra N$ satisfying $\delta_M(f\cdot \xi \cdot g) = \delta_{\cA}(f) \cdot \xi g + f \cdot \delta_M(\xi) \cdot g + f \xi \cdot \delta_{\cA}(g)$, for $f \in {}_{x_3}\cA_{x_2}$, $\xi \in {}_{x_2}M_{x_1}$, $g \in {}_{x_1}\cA_{x_0}$.
From the universal property of $\Omega^1$, this implies that
\begin{equation*}
	\Omega^1(\Pi) \, = \, \Big( \, ( \, \Pi \otimes_{\cA} \Omega^1(\cA) \otimes_{\cA} \Pi \,) \, \oplus \, ( \,\Pi \otimes_{\scO} DM \otimes_{\scO} \Pi  \,) \, \Big) \, / \, (D(f\cdot \xi \cdot g) = Df \cdot \xi g + f \cdot D\xi \cdot g + f \xi \cdot Dg)
\end{equation*}
where we have written $\Pi \otimes_{\scO} DM \otimes_{\scO} \Pi$ to mean the $\Pi$-bimodule $\Pi \otimes_{\scO} M \otimes_{\scO} \Pi$, but with elements $\xi \in {}_{y}M_{x}$ written as $D\xi \in {}_{y}DM_{x}$.

After setting the component $\Pi \otimes_{\cA} \Omega^1(\cA) \otimes_{\cA} \Pi$ to be zero, the remaining relation becomes $(D(f\cdot \xi \cdot g) = f \cdot D\xi \cdot g$, so that the quotient is precisely $\Pi \otimes_{\cA} DM \otimes_{\cA} \Pi = \iota^e_!(M)$. This shows \eqref{TAM_Omega_SES}.

If $M$ is cofibrant, then it is a retract of a cellular bimodule $P$. Thus, in particular $P$ is free as a graded bimodule over a set $\{ \eta_i \in {}_{y_i}P_{x_i} \}_{i \in I}$. We will now forget the differentials and regard $\cA$, $\Pi$ and $T_{\cA}P$ as graded $k$-categories. Then there is a retract $\Pi \xra{i} T_{\cA}P \xra{r} \Pi$, commuting with the canonical maps from $\cA$. Denote by $j : \cA \ra T_{\cA} P$ the canonical inclusion, then we have 
\begin{equation*}
	\Omega^1(T_{\cA} P) = j^e_!(\Omega^1(\cA)) \, \oplus \, \bigoplus_{i \in I} \cA_{y_i} \cdot D\eta_i \cdot  {}_{x_i}\cA
\end{equation*}
Applying $r^e_!$, we see that the canonical map $\iota^e_!(\Omega^1(\cA)) \ra r^e_!(\Omega^1(T_{\cA}P))$ is injective. Since the map $\iota^e_!(\Omega^1(\cA)) \ra \Omega^1(\Pi)$ is a retract of this map, it is also injective.
\epf


\bcor  \label{TAM_cofib_cor}
Assume that $M \in \Mod(\cAe)$ is cofibrant, then 
\begin{enumerate}
	\item The canonical maps in \eqref{TAM_S_SES} can be extended to an exact triangle in $\cD(\Pie)$:
	\begin{equation}   \label{TAM_S_triangle}
		\ldots \raq \iota^e_! \, \cS(\cA) \xraq{\gamma} \cS(\Pi) \xraq{\rho} \iota^e_! \, M[1]  \raq \ldots
	\end{equation}
    \item The canonical dg functor $\iota : \cA \rinto \Pi$ is a cofibration.
    \item If $\cA$ is almost cofibrant, then so is $\Pi$.
    \item If $\cA$ is smooth, and $M$ is perfect, then $\Pi$ is smooth.
\end{enumerate}
\ecor

\bpf

(1) follows from Proposition \ref{Omega_Pi_SES}.

For (2), notice that if $M$ is cellular, then $\iota : \cA \ra \Pi$ is a cellular extension. In general, since we assume that $M$ is cofibrant, it is a retract of a cellular $N$. Thus, $\iota : \cA \ra \Pi$ is a retract of a cellular extension $\iota' : \cA \ra \Pi'$ in the undercategory $\cA \downarrow \dgcat_k$, and is therefore a cofibration.

For (3), to show that $\cS(\Pi)$ is cofibrant, use the same argument as in (2) to argue that the canonical map $\iota^e_!(\cSA) \ra \cS(\Pi)$ is a cofibration. Namely, consider the bimodule retract \eqref{Omega_retract_functorial} induced for the retract $\Pi \ra \Pi' \ra \Pi$, then we see that $\iota^e_!(\cSA) \ra \cS(\Pi)$ is a retract of a cellular extension in the undercategory $\iota^e_!(\cSA) \downarrow \Mod(\Pie)$.

The fact that $\Pi$ is linearly cofibrant also follows from the same retract argument. Namely, a cellular extension $\cA \ra \Pi'$ clearly remains linearly cofibrant, and hence so is its retract $\Pi$.

For (4), simply consider the exact triangle \eqref{TAM_S_triangle} and notice that both $\iota^e_! \, \cS(\cA)$ and $\iota^e_! \, M$ are perfect.
\epf

In the above, we have defined the dg category $\Pi = (T_{\cA}M , d_{\eta})$ from the data $(M,\eta)$, where $M \in \Mod(\cAe)$ and $\eta : M[-1] \ra \cA$ is a map in $\Mod(\cAe)$. We now investigate the case when $[M] \in \cD(\cAe)$ and $\eta : [M][-1] \ra [\cA]$ is a map in $\cD(\cAe)$.

\bpp  \label{Pi_M_eta_homotopy}
Given $[M] \in \cD(\cAe)$ and a map $\eta : [M][-1] \ra [\cA]$ in $\cD(\cAe)$, choose a cofibrant $M \in \Mod(\cAe)$ and a map $\eta : M[-1] \ra \cA$ in $\Mod(\cAe)$ that represents $([M],[\eta])$. Then the quasi-isomorphism%
\footnote{We say that a dg functor $F : \cA \ra \cB$ is a \emph{quasi-isomorphism} if it is bijective on object sets, and if $F : \cA(x,y) \ra \cB(F(x),F(y))$ is a quasi-isomorphism for all $x,y \in \Ob(\cA)$.} type of the resulting dg category $\Pi := (T_{\cA}M , d_{\eta})$ is independent of the choice of the representative $(M,\eta)$.
\epp

\bpf
We first fix $M$ and show that the isomorphism type of $\Pi$ is independent of the choice of $\eta$. Thus, if $\eta' : M[-1] \ra \cA$ represents the same map in $\cD(\cAe)$, then since $M$ is cofibrant, $\eta$ and $\eta'$ are homotopic to each other. {\it i.e.,} there exists $h \in \Homcom_{1}(M[-1],\cA) = \Homcom_0(M,\cA)$ such that $d(h) = \eta - \eta'$. (Here, we regard $h$ as a degree $0$ element in $\Homcom(M,\cA)$ so that $d(h) = d \circ h - h \circ d$). There is a unique map $\varphi_h : T_{\cA}(M) \ra T_{\cA}(M)$ of graded $k$-categories which is the identity on $\cA$ and is given by $(h,\id) : M \ra \cA \oplus M \subset T_{\cA}(M)$ on $M$. Since $\varphi_h \circ \varphi_{-h} = \varphi_{-h} \circ \varphi_{h} = \id$, we see that $\varphi_h$ is an isomorphism of graded $k$-categories. Moreover, the assumption $dh - hd = \eta - \eta'$ implies that $\varphi_h \circ d_{\eta} = d_{\eta'} \circ \varphi_{h}$ (check on the generators $\cA$ and $M$). Thus, $\varphi_h : (T_{\cA}M, d_{\eta}) \xra{\cong} (T_{\cA}M, d_{\eta'})$ is an isomorphism of dg categories.

Next, we show that the quasi-isomorphism type of $\Pi$ is independent of the choice of cofibrant $M$. Given any other cofibrant representative $M'$, there is a quasi-isomorphism $f : M' \xra{\sim} M$. Since we have shown above that the isomorphism type of $\Pi$ is independent of the choice of the representative $\eta$, we may assume that $\eta' : M \ra \cA$ is given by the composition $M' \xra{f} M \xra{\eta} \cA$. Hence, the map of graded $k$-categories $\varphi : T_{\cA}M' \ra T_{\cA}M$ given by identity on $\cA$ and $f$ on $M'$ satisfies $\varphi \circ d_{\eta'} = d_{\eta} \circ \varphi$, so that it is a dg functor.
Recall that we have assumed that $\cA$ is $k$-flat. Thus, any cofibrant $M$ is in particular left $k$-flat, so that the map $M' \otimes_{\cA} \stackrel{(n)}{\ldots } \otimes_{\cA} M' \ra M \otimes_{\cA} \stackrel{(n)}{\ldots } \otimes_{\cA} M$ is a quasi-isomorphism for each $n \geq 0$. To show that $\varphi$ is a quasi-isomorphism, it suffices to consider the filtration $F_{\bullet}$ on the chain complex $T_{\cA}(M)(x,y)$ given by $F_r := \bigoplus_{n \leq r} \, ({}_y M \otimes_{\cA} \stackrel{(n)}{\ldots } \otimes_{\cA} M_{x})$, and similarly on $T_{\cA}(M')(x,y)$. The map $\varphi$ preserves this filtration, and the induced map on the associated graded complexes is precisely $\varphi : (T_{\cA}M', d_0) \ra (T_{\cA}M, d_0)$, which we have just shown to be a quasi-isomorphism. By a standard spectral sequence argument, we see that $\varphi : (T_{\cA}M', d_{\eta'}) \ra (T_{\cA}M, d_{\eta})$ is also a quasi-isomorphism.
\epf

\bdf
A \emph{tensor extension data} is a triple $(\cA,M,\eta)$ where $\cA \in \dgcatflat_k$ is a $k$-flat small dg category, $M \in \cD(\cAe)$, and $\eta: M[-1] \ra \cA$ is a map in $\cD(\cAe)$. The data $(M,\eta)$ is also identified as an object in the undercategory $\cD(\cAe) \downarrow \cA$.

A \emph{map of extension data}  $(F,f) : (\cA,M,\eta) \ra (\cB,N,\xi)$ consists of a dg functor $F : \cA \ra \cB$ together with a map $f_L : {\bm L}F^e_!(M) \ra N$ in $\cD(\cBe)$, or equivalently a map $f_R : M \ra F^{e*}(N)$ in $\cD(\cAe)$, such that one, and hence the other, of the following diagrams commute:
\begin{equation}  \label{fL_fR_comm_diag}
	\begin{tikzcd}
		{\bm L}F^e_!(M)[-1] \ar[r, "f_L"] \ar[rr, bend left, "\gamma_F \circ {\bm L}F^e_!(\eta)"] & N[-1] \ar[r, "\xi"] & \cB 
&
&M[-1] \ar[r, "f_R"] \ar[rr, bend left, "F \circ \eta"] & F^{e*}(N)[-1] \ar[r, "F^{e*}(\xi)"] & F^{e*}(\cB)
\end{tikzcd}
\end{equation} 
\edf

Proposition \ref{Pi_M_eta_homotopy} allows one to define the dg category $(T_{\cA}M , d_{\eta})$ out of a tensor extension data, well-defined up to quasi-isomorphism. This construction is functorial: a map $(F,f) : (\cA,M,\eta) \ra (\cB,N,\xi)$ of extension data induces a dg functor $\widetilde{F} : (T_{\cA} M, d_{\eta}) \ra (T_{\cB}N , d_{\xi})$. 
Indeed, notice that for $F = \id : \cB \ra \cB$, the functoriality in $(M,\eta)$ follows from Proposition \ref{Pi_M_eta_homotopy}. The case for general $F$ then follows by combining this with the ``initial case'', which we discuss now.

Given a tensor extension data  $(\cA,[M],[\eta])$ and a dg functor $F: \cA \ra \cB$ to another $k$-flat small dg category $\cB$. Then there is a pushforward tensor extension data on $\cB$, given by ${\bm L}F^e_!(M) \in \cD(\cBe)$, and $[\eta'] : {\bm L}F^e_!(M)[-1] \xra{{\bm L}F^e_!(\eta)} {\bm L}F^e_!(\cA) \xra{\gamma_F} \cB$ in $\cD(\cBe)$. 
If we choose cofibrant representative $M$ of $[M]$ and $\eta : M[-1] \ra \cA$ that represents $[\eta]$, then the data on $\cB$ may be represented by $F^e_!(M)$ and $\eta' : F^e_!(M)[-1] \xra{F^e_!(\eta)} F^e_!(\cA) \ra \cB$ in $\Mod(\cBe)$. 
Using these choices, there is an obvious dg functor
\begin{equation*}
	\widetilde{F} \, : \, \Pi \, := \, ( \, T_{\cA}(M)  \, ,\,  d_{\eta} \, ) 
	\raq ( \, T_{\cB}(F^e_!(M)) \, , \, d_{\eta'} \, ) \, =: \, \Pi' 
\end{equation*}

\bpp  \label{TAM_homotopy_invariance}
\begin{enumerate}
	\item If $F$ is quasi-fully faithful, then so is $\widetilde{F}$.
	\item If $F$ is a quasi-equivalence, then so is $\widetilde{F}$.
	\item If $F$ is a derived equivalence, then so is $\widetilde{F}$.
\end{enumerate}
\epp

\bpf
Writing $M' = F^e_!(M)$, then we see that
\begin{equation*}
	M' \otimes_{\cB} \ldots \otimes_{\cB} M' \, = \, \cB \otimes_{\cA} M  \otimes_{\cA} \cB  \otimes_{\cA} \ldots  \otimes_{\cA} \cB  \otimes_{\cA} M  \otimes_{\cA} \cB
	\, = \, \cB^{\otimes n+1} \otimes_{(\cAe)^{\otimes n}} M^{\otimes n}
\end{equation*}
When written in this form, the map $\widetilde{F} : \Pi(x,y) \ra \Pi'(F(x),F(y))$ is given on the $n$-th component by
\begin{equation}  \label{F_tilde_nth_comp}
	({}_y\cA \otimes \cA^{\otimes n-1} \otimes \cA_x ) \otimes_{(\cAe)^{\otimes n}} M^{\otimes n} \raq ({}_{F(y)}\cB \otimes \cB^{\otimes n-1} \otimes \cB_{F(x)} ) \otimes_{(\cAe)^{\otimes n}} M^{\otimes n}
\end{equation}

For (1), to show that $\widetilde{F}$ is quasi-fully faithful, it suffices, by a spectral sequence argument as in the proof of Proposition \ref{Pi_M_eta_homotopy}, to assume that $\eta = \eta' = 0$ so that it suffices to show that \eqref{F_tilde_nth_comp} is a quasi-isomorphism. Since $F$ is assumed to be quasi-fully faithful, and $\cA$ and $\cB$ assumed to be $k$-flat, the map ${}_y\cA \otimes \cA^{\otimes n-1} \otimes \cA_x \ra {}_{F(y)}\cB \otimes \cB^{\otimes n-1} \otimes \cB_{F(x)}$ is a quasi-isomorphism of $(\cAe)^{\otimes n}$-modules. 
Since $M$ is a cofibrant bimodule, $M^{\otimes n}$ is a cofibrant $(\cAe)^{\otimes n}$-module, therefore \eqref{F_tilde_nth_comp} is indeed a quasi-isomorphism.

For (2) and (3), we have already seen that $\widetilde{F}$ is quasi-fully faithful in (1). To verify the remaining conditions, we consider the following two commutative (up to natural isomorphisms) diagrams of functors:
\begin{equation*}
	\begin{tikzcd}
		H_0(\cA) \ar[r] \ar[d] & H_0(\cB) \ar[d] & & \Dperf(\cA) \ar[r] \ar[d] & \Dperf(\cB) \ar[d] \\
		H_0(\Pi) \ar[r] & H_0(\Pi') & & \Dperf(\Pi) \ar[r] & \Dperf(\Pi')
	\end{tikzcd}
\end{equation*}

For (2), use the diagram on the left. By assumption $H_0(\cA) \ra H_0(\cB)$ is essentially surjective. Since $H_0(\cB) \ra H_0(\Pi')$ is clearly essentially surjective (in fact, it is bijective on object sets), so is their composition $H_0(\cA) \ra H_0(\Pi')$. In particular, this implies that $H_0(\Pi) \ra H_0(\Pi')$.

For (3), we use the same argument, except that the notion of essential surjectivity on $H_0$ is replaced by the following notion: we say that a dg functor $F: \cA_1 \ra \cA_2$ has dense image if the smallest split-closed triangulated subcategory of $\Dperf(\cA_2)$ containing the essential image of ${\bm L}F_! : \Dperf(\cA_1) \ra \Dperf(\cA_2)$ is $\Dperf(\cA_2)$.
Notice that: (i) if $F: \cA_1 \ra \cA_2$ and $G: \cA_2 \ra \cA_3$ have dense images, then so does $G \circ F$; (ii) if $F : H_0(\cA_1) \ra H_0(\cA_2)$ is essentially surjective, then $F$ has dense image; (iii) if $G \circ F$ has dense image, then so does $G$. 
Using these three properties, one can use the same argument as in (2) to establish (3).
\epf

\bdf
Let $\cA \in \dgcatflat_k$ be a small $k$-flat dg category. Assume that $\cA$ is smooth. Given a negative cyclic class $\widetilde{\eta} \in \HN_{n-2}(\cA)$, with underlying Hochschild class $\eta \in \HH_{n-2}(\cA)$, so that $\eta$ determines a map $\eta^{\#} : \cA^![n-2] \ra \cA$ in $\cD(\cAe)$. 
Then its \emph{deformed $n$-Calabi-Yau completion} with deformation parameter $\eta$ is the dg category $\Pi_n(\cA,\eta)$ defined as in Proposition \ref{Pi_M_eta_homotopy} from the data $M = \cA^![n-1]$ and the map $\eta^{\#}$.
With a slight abuse of notation, we may write
\begin{equation*}
	\Pi_n(\cA,\eta) \, = \, ( \, T_{\cA}(\cA^![n-1]) \, , \, d_{\eta} \, ) 
\end{equation*}
Notice that the deformed $n$-Calabi-Yau completion depends only on $\eta$, but we always require it to have a given negative cyclic lift $\widetilde{\eta}$. See Remark \ref{CYcom_not_CY_remark} below for a reason.
\edf

Given a dg functor $F : \cA \ra \cB$ between small $k$-flat dg categories. For a given $\eta \in \HH_{n-2}(\cA)$, denote by $\eta'' \in HH_{n-2}(\cB)$ its image under $F$. In fact, the canonical map $\HH_{n-2}(\cA) \ra \HH_{n-2}(\cB)$ factors as follows:
\begin{equation*}
	\cA \otimes^{{\bm L}}_{\cAe} \cA \raq \cA \otimes^{{\bm L}}_{\cAe} \cB \, \simeq \, {\bm L}F^e_!(\cA) \otimes^{{\bm L}}_{\cBe} \cB \raq \cB \otimes^{{\bm L}}_{\cBe} \cB
\end{equation*}
Let $\eta' \in H_{n-2}( {\bm L}F^e_!(\cA) \otimes^{{\bm L}}_{\cBe} \cB )$ be the image under $\eta$ in the first map. 
Then the three elements $\eta$, $\eta'$, $\eta''$ induce the three downward pointing maps in the following commutative diagram in $\cD(\cB)$:
\begin{equation}  \label{CY_com_compare_1}
	\begin{tikzcd}
		{\bm L}F^e_!(\cA^!)[n-2]   \ar[rd, "\gamma_F \circ {\bm L}F^e_!(\eta^{\#})"'] \ar[r, "\eqref{F_shriek_dual_1}"]
		&  {\bm L}F^e_!(\cA)^![n-2]  \ar[d, "(\eta')^{\#}_L"]
		& \cB^! \ar[ld, "(\eta'')^{\#}"] \ar[l, "(\gamma_F)^!"'] \\
		& \cB
	\end{tikzcd}
\end{equation}

Now if $\cA$ is smooth, then the first horizontal map is an isomorphism. If $F$ is a derived equivalence, then the second horizontal map is an isomorphism (see Proposition \ref{derived_equiv_criteria}). Combining with Propositions \ref{Pi_M_eta_homotopy} and \ref{TAM_homotopy_invariance}, we have

\bpp  \label{CY_com_homotopy_inv}
The quasi-isomorphism/quasi-equivalence/derived equivalence type of the deformed $n$-Calabi-Yau completion $\Pi_n(\cA,\eta)$ only depends on the quasi-isomorphism/quasi-equivalence/derived equivalence type of the data $(\cA,\eta)$ where $\cA \in \dgcatflat_k$ is smooth and $\eta \in \HH_{n-2}(\cA)$.  
\epp

\bpf
Since the horizontal maps in \eqref{CY_com_compare_1} are isomorphisms, we see from Propositions \ref{Pi_M_eta_homotopy} that the dg categories $(T_{\cB}(N) , d_{\xi})$ are quasi-isomorphic when $(N,\xi)$ is taken to be either of the three downward arrows in \eqref{CY_com_compare_1}. Then Proposition \ref{TAM_homotopy_invariance} says that $(T_{\cA}(\cA^![n-1]), d_{\eta})$ and $(T_{\cB}(N) , d_{\xi})$ are quasi-isomorphic/quasi-equivalent/derived equivalent for $(N,\xi)$ being the leftmost downward arrow in \eqref{CY_com_compare_1}.
\epf

When working with explicit examples, it is convenient to assume that $\cA$ is almost cofibrant (by Proposition \ref{CY_com_homotopy_inv} we may replace $\cA$ by a cofibrant resolution). Then $\cA^!$ can be explicitly represented by $\cSA^{\vee} \in \Mod(\cAe)$. It is also convenient to represent the deformation parameter by cycles in the X-complex $\eta \in Z_{n-2}(X(\cA))$. 
Thus, if we choose any cofibrant resolution $\cP \xra{\sim} \cSA^{\vee}$ in $\Mod(\cAe)$, then this gives an explicit model $(T_{\cA}(\cP[n-1]) , d_{\eta})$ for the deformed $n$-Calabi-Yau completion.
In particular, when $\cA$ is finitely cellular, then we may take $\cP = \cSA^{\vee}$. In this case, the explicit model $(T_{\cA}(\cSA^{\vee}[n-1]) , d_{\eta})$ remains finitely cellular.

As standard examples, recall that both the (derived) deformed preprojective algebras and the Ginzburg dg algebras can be realized as deformed Calabi-Yau completions by using this explicit model $(T_{\cA}(\cSA^{\vee}[n-1]) , d_{\eta})$.

\beg  \label{2CY_com_eg1}
Let $A = k[x]$. Then the short resolution $\cS(A)$ is given by
\begin{equation*}
\cS(A) \, = \, [ \, A \cdot Dx \cdot A \xraq{\alpha} A \cdot E \cdot A ]
\end{equation*}
where we have $\alpha(f(x) \cdot Dx \cdot g(x)) = f(x)x \cdot E \cdot g(x) - f(x) \cdot E \cdot xg(x)$.
In particular, the X-complex of $A$ is given by
\begin{equation*}
X(A) \, = \, [ \,  A \cdot Dx \xraq{0} A \cdot E \, ]
\end{equation*}

Thus, any polynomial $f(x) \in A$ represents a class $\eta = f(x) \cdot E \in \HH_0(A)$. It can be shown that $\eta$ has a negative cyclic lift if and only if $Df \in A \cdot Dx$ is zero. For example, if $\bQ \subset k$, then this is true if and only if $f(x) \in k \subset k[x]$.

Using the explicit model $\cS(A)^{\vee}$ for $A^!$, we obtain an explicit model for the ``deformed $2$-Calabi-Yau completion'' (we add quotation marks because we carry out this construction for $\eta$ that do not necessarily have a negative cyclic lift)
\begin{equation}  \label{2CY_com_alg_1}
\Pi_2(A,\eta) \, = \, ( \, k\langle x,y,t \rangle \, , \quad  d(t) = xy-yx+f(x) \, )
\end{equation}
where $|x|=|y|=0$ and $|t|=1$.
Namely, we rename $(Dx)^{\vee}$ into $y$, and $E^{\vee}[1]$ into $t$. 
Notice that the differential in $\cS(A)^{\vee}$ is given by $d(E^{\vee}) = x \cdot (Dx)^{\vee} - (Dx)^{\vee} \cdot x$ (see \eqref{dualmod_diff}), which gives rise to the part $xy-yx$ in the differential $d(t)$, the remaining part $f(x)$ comes from the deformation parameter $\eta = f(x)\cdot E$.

\eeg

We now consider the relative case.

Given a dg functor $F: \cA \ra \cB$ between small $k$-flat dg categories, and given a class $\eta \in \HH_{n-2}(\cB,\cA)$, then there is an induced map of exact triangles \eqref{rel_Hoch_induce_map_1}. We will focus on the square on the right hand side. Assume that $\cA$ is smooth, then by the discussion that follows Definition \ref{rel_CY_def}, we may replace ${\bm L}F^e_!(\cA)^!$ by ${\bm L}F^e_!(\cA^!)$, and the square becomes
\begin{equation}  \label{rel_Hoch_compatible}
\begin{tikzcd}
{\bm L}F^e_!(\cA^!)[n-3] \ar[r] \ar[d, "{\bm L}F^e_!((\eta_{\cA})^{\#})"'] & 
\cone [ \cB^! \ra {\bm L}F^e_!(\cA)^! ][n-3]  \ar[d, "(\eta''')^{\#}_L"] \\
{\bm L}F^e_!(\cA) \ar[r, "\gamma_F"]
& \cB
\end{tikzcd}
\end{equation}
where $\eta_{\cA} \in \HH_{n-3}(\cA)$ is the image of $\eta$ under the connecting homomorphism $\HH_{n-2}(\cB,\cA) \ra \HH_{n-3}(\cA)$.

Notice that \eqref{rel_Hoch_compatible} is precisely the compatibility condition \eqref{fL_fR_comm_diag} for the map of extension data from $(\cA,\cA^![n-2],\eta_{\cA}^{\#})$ to $(\cB,N,\xi)$, where $N = \cone [ \cB^! \ra {\bm L}F^e_!(\cA)^! ][n-2]$ and $\xi = (\eta''')^{\#}_L$.
Therefore, there is a dg functor
\begin{equation}   \label{rel_CY_com_map}
\Pi_{n-1}(\cA, \eta_{\cA}) \raq \Pi_n(\cB,\cA, \eta)\,  , \text{ where } 
\begin{cases*}
\Pi_{n-1}(\cA, \eta_{\cA}) \, := \, (\, T_{\cA}(\cA^![n-2]) \, , \, d_{\eta_{\cA}} \,) \\
\Pi_n(\cB,\cA, \eta) \, := \, (\, T_{\cB}( \cone [ \cB^! \ra {\bm L}F^e_!(\cA)^! ][n-2] ) \, , \, d_{\eta'''} \,)
\end{cases*} 
\end{equation}

\bdf  \label{relCYcom_def}
Given a dg functor $F : \cA \ra \cB$ between small $k$-flat dg categories, both of which are smooth, and given a relative negative cyclic class $\widetilde{\eta} \in \HN_{n-2}(\cB,\cA)$ with underlying relative Hochschild class $\eta \in \HH_{n-2}(\cB,\cA)$, then the induced dg functor \eqref{rel_CY_com_map} is called the \emph{deformed relative $n$-Calabi-Yau completion} with deformation parameter $\eta$.

As in the absolute case, the deformed relative $n$-Calabi-Yau completion depends only on the relative Hochschild class $\eta \in \HH_{n-2}(\cB,\cA)$, but we always require that it has a given negative cyclic lift.
\edf

As in the absolute case, this construction is invariant under quasi-isomorphism/quasi-equivalence/derived equivalence. The detailed formulation will be omitted.

\brm  \label{CYcom_not_CY_remark}
Given a small $k$-flat dg category $\cA$, then we may say that $\cA$ is weakly $n$-Calabi-Yau if there is an isomorphism $\cA^![n] \cong \cA$ in $\cD(\cAe)$. 
Example \ref{2CY_com_eg1} gives an example that, when the deformation parameter does not have a negative lift, then the ``deformed $n$-Calabi-Yau completion'' $\Pi_n(\cA,\eta)$ may not be weakly $n$-Calabi-Yau (here we consider $n=2$).

Namely, take $f(x) = x$, and assume that $\bQ \subset k$. Denote by $\Pi$ the dg algebra \eqref{2CY_com_alg_1}. We claim that $H_0(\Pi \otimes^{{\bm L}}_{\Pie} \Pi) \neq 0$ but $H_0(\Pi^![2] \otimes^{{\bm L}}_{\Pie} \Pi) = 0$, so that $\Pi$ cannot be weakly $2$-Calabi-Yau. (The reader is cautioned that in the previous version of the present paper, the calculation of $H_0(\Pi \otimes^{{\bm L}}_{\Pie} \Pi)$ and $H_0(\Pi^![2] \otimes^{{\bm L}}_{\Pie} \Pi)$ contained some mistakes.)

It is easy to see that $H_0(\Pi \otimes^{{\bm L}}_{\Pie} \Pi) \neq 0$. Namely, write $B := H_0(\Pi) = k\langle x,y\rangle / (xy-yx+x)$, then since $\Pi$ is non-negatively graded, we have 
$H_0(\Pi \otimes^{{\bm L}}_{\Pie} \Pi) = \HH_0(\Pi) = \HH_0(B) = B/[B,B]$.
There is a surjective map $B/[B,B] \ronto B/([B,B])$ to the commutativization of $B$. Since $B/([B,B]) = k[x,y]/(xy-yx+x) = k[x,y]/(x) = k[y]$, we see that $B/[B,B]$ is non-zero. 
In fact, by \eqref{B_xy_commutator} it can be shown that, if $\bQ \subset k$ then the map $B/[B,B] \ra B/([B,B]) = k[y]$ is an isomorphism of $k$-modules.

The short resolution $\cS(\Pi)$ of $\Pi$ has a basis given by
$\{ E , \, sDx, \, sDy, \, sDt  \}$, with differentials
\begin{equation*}
\begin{split}
d(E) &= 0 \qquad \quad
d(sDx) = x \cdot E - E \cdot x  \qquad \quad 
d(sDy) = y \cdot E - E \cdot y \\
d(sDt) &= t \cdot E - E \cdot t - sD(xy - yx + x) \\
&= t \cdot E - E \cdot t - sDx \cdot y - x \cdot sDy + sDy \cdot x + y \cdot sDx - sDx
\end{split}
\end{equation*}

By taking the dual, we see that $\cS(\Pi)^{\vee}$ 
has basis 
$\{ E^{\vee} , \, (sDx)^{\vee}, \, (sDy)^{\vee}, \, (sDt)^{\vee}  \}$,
with differentials (see \eqref{dualmod_diff}) given by
\begin{equation*}
\begin{split}
d(sDt)^{\vee} &= 0 \qquad \quad
d(sDy)^{\vee} = - x \cdot (sDt)^{\vee} + (sDt)^{\vee} \cdot x \\
d(sDx)^{\vee} &= y \cdot (sDt)^{\vee} - (sDt)^{\vee} \cdot y + (sDt)^{\vee} \\
d(E^{\vee}) &= -(sDx)^{\vee} \cdot x + x \cdot (sDx)^{\vee} -(sDy)^{\vee} \cdot y + y \cdot (sDy)^{\vee} 
-(sDt)^{\vee} \cdot t + t \cdot (sDt)^{\vee}
\end{split}
\end{equation*}

In particular, the complex $\cS(\Pi)^{\vee}[2] \otimes_{\Pie} \Pi$ is given by
\begin{equation*}
\begin{tikzcd}[row sep = 10]
& \Pi \cdot s^2(sDy)^{\vee}  \ar[dr, "\varphi_1"] & \\
\Pi \cdot s^2 E^{\vee} \ar[ur] \ar[dr]  \ar[rr]
& & \Pi \cdot s^2(sDt)^{\vee} \\
& \Pi \cdot s^2(sDx)^{\vee} \ar[ur, "\varphi_2"']
\end{tikzcd}
\end{equation*}
where, in particular, the maps $\varphi_1$ and $\varphi_2$ are given by 
\begin{equation*}
\begin{split}
\varphi_1 ( g \cdot s^2 (sDy)^{\vee} ) \, &= \, (- g \cdot x + x \cdot g) s^2(sDt)^{\vee} \\
\varphi_2 ( g \cdot s^2 (sDx)^{\vee} ) \, &= \, (g \cdot y - y \cdot g + g) s^2(sDt)^{\vee}
\end{split}
\end{equation*}

Again, write $B := H_0(\Pi) = k\langle x,y\rangle / (xy-yx+x)$, then we have
\begin{equation}  \label{B_quot_g12}
H_0( \cS(\Pi)^{\vee}[2] \otimes_{\Pie} \Pi )
\, = \, B / \{ - g_1 \cdot x + x \cdot g_1 + g_2 \cdot y - y \cdot g_2 + g_2 \}_{g_1,g_2 \in B}
\end{equation}

Notice that $B$ is the universal enveloping algebra of the Lie algebra $(kx \oplus ky, [x,y]=x)$, so that $\{x^n y^m\}_{m,n \in \bN}$ is a $k$-linear basis of $B$.
One can compute the terms in the quotient \eqref{B_quot_g12} by taking $g_1$ and $g_2$ to be these monomials. Notice that in $B$ we have
\begin{equation}  \label{B_xy_commutator}
[y, x^n y^m] \, = \, nx^n y^m  \qquad \text{ and } \qquad [x,x^n y^m]  \, \equiv \, mx^{n+1}y^{m-1} \quad \text{ mod } \oplus_{i=0}^{m-2} k x^{n+1}y^i
\end{equation}
In particular, since we assume that $\bQ \subset k$, the terms $-g_1 \cdot x + x \cdot g_1$ span every monomial $x^ny^m$ with $n > 0$. Moreover, the terms $g_2 \cdot y - y \cdot g_2 + g_2$ span every monomial $y^m$ by taking $g_2 = y^m$. Together, this implies that \eqref{B_quot_g12} is zero. 
\erm

\section{Main theorem: the absolute case}  \label{sec_main_thm_abs}

\bthm  \label{main_thm_abs}
For any smooth $k$-flat small dg category $\cA$ and any negative cyclic class $\widetilde{\eta} \in \HN_{n-2(\cA)}$, the deformed $n$-Calabi-Yau completion $\Pi_n(\cA,\eta)$ has an $n$-Calabi-Yau structure%
\footnote{See also Remark \ref{CY_str_can_remark} for a discussion about the extent to which this $n$-Calabi-Yau structure is canonical.}.

Moreover, if the deformation parameter $\widetilde{\eta}$ is exact in the sense that it is in the image of $B : \HC_{n-3}(\cA) \ra \HN_{n-2}(\cA)$, then this $n$-Calabi-Yau structure on $\Pi_n(\cA,\eta)$ is exact.
\ethm

We will first prove this theorem for the case when $\cA$ is finitely cellular and the deformation parameter $\widetilde{\eta}$ is zero. Then we will briefly discuss what needs to be modified in order to cover the case when $\cA$ is finitely cellular and the the deformation parameter $\widetilde{\eta}$ is non-zero. Then we will proceed to prove the general case, which essentially follows the same argument as for the finitely cellular case.
The readers only interested the general statement may skip directly to its proof as it is logically independent of the finitely cellular case.
However, we have kept the proof for the finitely cellular case because in this case, the main ingredients of the proof admits descriptions in terms of explicit formulas (see \eqref{J_theta_dual_basis} and \eqref{B_J_theta_dual_basis}), so that it showcases the formulaic nature of the construction of Calabi-Yau completions (which we think is a major merit of this construction), and at the same time aid the reader's understanding of the proof in the general case.

As a preparation for the proof of Theorem \ref{main_thm_abs}, we give the following lemma about some properties of ``Casimir elements''. Strictly speaking, this lemma is not needed in the finitely cellular case, but we still present it first for the sake of logical clarity.

\blm  \label{Casimir_element_lemma}
Let $\cA$ be a $k$-flat small dg category. Let $Q \in \Mod(\cAe)$ be perfect and cofibrant, and let $p : P \xra{\sim} Q^{\vee}$ be a cofibrant replacement. Denote by $p^a$ the map $p^a : Q \xra{u} Q^{\vee \vee} \xra{p^{\vee}} P^{\vee}$. Then
\begin{enumerate}
	\item There exists an element $\theta \in Z_0(P \otimes_{\cAe} Q)$ such that the induced maps
	\begin{equation*}
		\begin{split}
			\varphi' \ &: \, P \xraq{p} Q^{\vee} \xraq{\theta^{\#}_R} P \\
			\varphi'' \ &: \, Q \xraq{p^a} P^{\vee} \xraq{\theta^{\#}_L} Q
		\end{split}
	\end{equation*}
	are both homotopic to the identities. In fact, $\varphi' \simeq \id$ if and only if $\varphi'' \simeq \id$.
	Any such element $\theta$ is called a Casimir element for $(Q,P,p)$.
	\item Let $\theta \in Z_0(P \otimes_{\cAe} Q)$ be a Casimir element, then for any closed element $\xi \in Z_m(Q_{\natural})$, denote by $\xi'$ the image of $\theta$ under the map $P \otimes_{\cAe} Q \xra{p \otimes \id} Q^{\vee} \otimes_{\cAe} Q \xra{\xi^{\#} \otimes \id} \cA[-m] \otimes_{\cAe} Q = Q_{\natural}[-m]$. Then $\xi'$ is homologous to $\xi$. 
	\item If $Q$ is strictly perfect (hence in particular cofibrant), then we may take $p = \id$. In this case, we can require $\varphi' = \id$ and $\varphi'' = \id$ ({\it i.e.,} they are not just homotopic, but equal). In this case, we call $\theta$ a strict Casimir element. For a strict Casimir element, we have strict equality $\xi' = \xi$. 
	\item If $F : \cA \ra \cB$ is a dg functor to another $k$-flat small dg category. Then the map $\widetilde{p} : F^e_!(P) \ra F^e_!(Q)^{\vee}$ defined as in \eqref{F_shriek_dual_2} is a cofibrant replacement, so that we may discuss about Casimir elements for $(F^e_!(Q),F^e_!(P),\widetilde{p})$. Let $j : P \otimes_{\cAe} Q \ra F^e_!(P) \otimes_{\Pie} F^e_!(Q)$ be the canonical map, then for any Casimir element $\theta$ for $(Q,P,p)$, $j(\theta)$ is a Casimir element for  $(F^e_!(Q),F^e_!(P),\widetilde{p})$. 
\end{enumerate}
\elm

\bpf
Since $Q$ is perfect, for any $M \in \Mod(\cAe)$, there is a canonical quasi-isomorphism
\begin{equation}  \label{M_tensor_Q_to_Hom}
	M \otimes_{\cAe} Q \raq \Homcom_{\cAe}(P,M) \, , \qquad \xi  \,\mapsto \, [P \xra{p} Q^{\vee} \xra{\xi^{\#}_R} M]
\end{equation}
since it is an explicit representative of the isomorphism $M \otimes^{{\bm L}}_{\cAe} Q \cong \RHomcom_{\cAe}(Q^!,M)$ in $\cD(k)$.

In particular, taking $M = P$, we can choose $\theta$ so that $\varphi'$ is homotopic to identity. Thus, to prove (1), it remains to show that $\varphi' \simeq \id$ if and only if $\varphi'' \simeq \id$. 
To show the implication ``$\Rightarrow$'', we notice that the following diagram is commutative:
\begin{equation}  \label{PQ_theta_dual}
	\begin{tikzcd}
		P^{\vee} \ar[r, "(\theta^{\#}_R)^{\vee}"] \ar[rd, "\theta^{\#}_L"']  & Q^{\vee \vee}  \ar[r, "p^{\vee}"] &  P^{\vee}  \\
		& Q \ar[u, "u"] \ar[ru, "p^a"']
	\end{tikzcd}
\end{equation}

Assume that $\varphi' \simeq \id$. Since the composition $p^{\vee} \circ (\theta^{\#}_R)^{\vee}$ is the dual of $\varphi' = \theta^{\#}_R \circ p$, it is homotopic to the identity. This shows that the composition $p^a \circ \theta^{\#}_L$ is the identity in the derived category $\cD(\cAe)$. Since $Q$ is perfect, $p^a$ is a quasi-isomorphism. By the uniqueness of inverse to $p^a$ in the derived category $\cD(\cAe)$, we see that  $\varphi'' = \theta^{\#}_L \circ p^a$ is also the identity in $\cD(\cAe)$. Since $Q$ is cofibrant, this is the same as saying that $\varphi''$ is homotopic to the identity.
This completes the proof that $\varphi' \simeq \id \, \Rightarrow \, \varphi'' \simeq \id$. The reverse implication ``$\Leftarrow$'' has a completely parallel proof: simply notice that $p^{aa} = p$.
This completes the proof of (1).

For (2), to show that $\xi'$ is homologous to $\xi$, it suffices to show that their images under \eqref{M_tensor_Q_to_Hom} for $M = \cA[-m]$ are homotopic. The image of $\xi$ under \eqref{M_tensor_Q_to_Hom} is given by $P \xra{p} Q^{\vee} \xra{\xi^{\#}} \cA[-m]$. The map \eqref{M_tensor_Q_to_Hom} is natural in $M$, so that the image of $\xi'$ under \eqref{M_tensor_Q_to_Hom} is given by $P \xra{\varphi'} P \xra{p} Q^{\vee} \xra{\xi^{\#}} \cA[-m]$. Since $\varphi'$ is homotopic to the identity, these two maps are homotopic.

For (3), since $Q$ is strictly perfect and $p = \id$, we see that \eqref{M_tensor_Q_to_Hom} is an isomorphism (instead of a quasi-isomorphism). We can then repeat the proof, but with homotopy replaced by equality.

For (4), simply notice that the composition
\begin{equation*}
	F^e_!(P) \xraq{\tilde{p}} F^e_!(Q)^{\vee} \xraq{ j(\theta)^{\#}_R } F^e_!(P)
\end{equation*}
is the map $F^e_!(\varphi')$, which is therefore homotopic to the identity.
\epf

\bpf[Proof of Theorem \ref{main_thm_abs} (finitely cellular, undeformed case)]
The proof in the finitely cellular, undeformed case can be motivated by the construction of the canonical symplectic structure on a cotangent bundle. Thus, let $M$ be a smooth manifold, with local coordinates $(q_1,\ldots,q_r)$, so that $T^*M$ has local coordinates $(q_1,\ldots,q_r,p_1,\ldots,p_r)$. One then verifies that the $1$-form $\lambda = \sum_{i=1}^r \, p_i \, dq_i$ is well-defined (independent of choice of coordiantes), and that $\omega = d\lambda$ is non-degenerate.
We now carry out a formal noncommutative analogue of this construction.

Assume that $\cA$ is finitely cellular, say $\cA = T_{\scO}(f_1,\ldots,f_r)$ with $f_i \in \cA(x_i,y_i)$, $|f_i| = m_i$, and with $d(f_i) \in T_{\scO}(f_1,\ldots,f_{i-1})$ for $i = 1, \ldots, r$. Then $\cSA$ is also finitely cellular, with a basis $\{ sDf_1 , \ldots, sDf_r\} \cup \{ E_x \}_{x \in \scO}$. 
Its bimodule dual is therefore also finitely cellular, with a basis $\{ (sDf_1)^{\vee} , \ldots, (sDf_r)^{\vee}\} \cup \{ E_x^{\vee} \}_{x \in \scO}$. 

Since $\cSA$ is finitely cellular, by Lemma \ref{Casimir_element_lemma}(3), there is a (unique) strictly Casimir element $\theta \in Z_0(\cSA^{\vee} \otimes_{\cAe} \cSA)$, characterized by the criterion that $\theta^{\#}_R : \cSA^{\vee} \ra \cSA^{\vee}$ is the identity (see Lemma \ref{Casimir_element_lemma}(3)). Clearly, we have
\begin{equation*}
	\theta \, = \, \Bigl( \sum_{i = 1}^r (-1)^{m_i+1} (sDf_i)^{\vee} \otimes (sDf_i) \Bigr) + \Bigl( \sum_{x \in \scO} E_x^{\vee} \otimes E_x \Bigr)
\end{equation*}
Notice that $\theta$ is closed because $\theta^{\#}_R$, a priori only a pre-map, is the identity, which is a map. 

Let $\Pi := T_{\cA}(\cSA^{\vee}[n-1])$, and let $\iota : \cA \ra \Pi$ be the canonical inclusion.
Consider the map 
\begin{equation}  \label{J_map_fin_cell}
	J \, : \, \cSA^{\vee}[n-1] \otimes_{\cAe} \cSA \rinto T_{\cA}(\cSA^{\vee}[n-1]) \otimes_{\cAe} \cSA  = 
	\Pi \otimes_{\Pie} (\iota^e_! (\cSA)) \xra{\gamma_{\natural}} \cS(\Pi)_{\natural}
\end{equation}

Inside $\Pi$, we will rewrite the generating arrows $s^{n-1}(sDf_i)^{\vee}$ by $g_i$. Thus $g_i \in \Pi(y_i,x_i)$ and $|g_i| = n-2-m_i$. We will also rewrite $s^{n-1}E_x^{\vee}$ by $c_x$. Thus $c_x \in \Pi(x,x)$ and $|c_x| = n-1$. We may write $\Pi = T_{\scO}(\{f_1,\dots,f_r,g_1,\dots, g_r\} \cup \{ c_x \}_{x \in \scO} )$.

Consider the image of $s^{n-1}\theta \in \cSA^{\vee}[n-1] \otimes_{\cAe} \cSA$ under $J$. By definition, it is given as
\begin{equation}  \label{J_theta_dual_basis}
	J(s^{n-1}\theta) \, = \, \Bigl( \sum_{i = 1}^r (-1)^{m_i+1} g_i \cdot (sDf_i) \Bigr) + \Bigl( \sum_{x \in \scO} c_x \cdot E_x \Bigr)
\end{equation}
This may be regarded as the formal noncommutative analogue of the $1$-form $\lambda = \sum_{i=1}^r \, p_i \, dq_i$ we discussed above (note the similar formula).

In the commutative case, the symplectic form is given by $d\lambda$, so that it depends on the de Rham differential $d : \Omega^1(M) \ra \Omega^2(M)$. We think of the mixed structure $B: X^{(1)}(\Pi) \ra X^{(2)}(\Pi)$ in \eqref{X_mixed_complex} as a noncommutative analogue of the de Rham differential. This is described in Lemma \ref{B_on_X1}. In particular, the element $\omega := B( J(s^{n-1}\theta) )$ is given by
\begin{equation}  \label{B_J_theta_dual_basis}
	B( J(s^{n-1}\theta) ) \, = \, \Bigl( \sum_{i = 1}^r (-1)^{m_i+1} sDg_i \otimes sDf_i + (-1)^{n(m_i+1)} sDf_i \otimes sDg_i \Bigr) + \Bigl( \sum_{x \in \scO} sDc_x \otimes E_x + E_x \otimes sDc_x\Bigr)
\end{equation}
By Theorem \ref{X_complex_thm_2} and \ref{X_complex_thm_3}, the element $\omega = B( J(s^{n-1}\theta))$ represents an exact negative cyclic class $[\omega] \in \HN_{n}(\Pi)$. It is clearly non-degenerate because the map $\omega^{\#} : \cS(\Pi)^{\vee}[n] \ra \cS(\Pi)$ is an isomorphism, as it sends a basis to a basis.
\epf

Before proceeding to the general case, we briefly discuss the finitely cellular, deformed case, in order to see what needs to be modified. Detailed will be subsumed in the proof of the general case below. 

Thus, let $\cA = T_{\scO}(f_1,\ldots,f_r)$ as above. By Theorem \ref{X_complex_thm_2}, we may represent the deformation parameter $\widetilde{\eta} \in \HN_{n-2}(\cA)$ by a closed element of degree $n-4$ in $F^1X^{\tot}(\cA)$. In other words, we may write $\widetilde{\eta}  =  \sum_{i=1}^{\infty}  \, \eta_i \cdot u^i$ for elements $\eta_i \in X^{(i)}(\cA)$ (see \eqref{eta_tilde_sum} below). 
The lowest order term $\eta_1$ will induce a differential in the deformed $n$-Calabi-Yau completion:
\begin{equation*}
	\Pi \, := \, ( \, T_{\cA}(\cSA^{\vee}[n-1]) \, , \, d_{\eta} \, )
\end{equation*}

Since the differential of $\Pi$ is deformed by $(\eta_1)^{\#}$, the pre-map \eqref{J_map_fin_cell}, while still well-defined, no longer commutes with differentials, since it involves the identitfication $\Pi = T_{\cA}(\cSA^{\vee}[n-1])$ which does not intertwine with the differentials.
In particular, $J(s^{n-1}\theta)$ is no longer $b$-closed, and is instead given by (see Lemma \ref{Casimir_element_lemma}(3))
\begin{equation}  \label{d_J_theta_fin_cell}
	b(J(s^{n-1}\theta)) \, = \, I(\eta_1)
\end{equation}
As such, $B( J(s^{n-1}\theta) ) \in X^{(2)}(\Pi)$ is also no longer $b$-closed. Instead, we have $b(B( J(s^{n-1}\theta) )) = -B(I(\eta_1)) = - I(B(\eta_1))$. At this point, we see why it is important to require that $\eta_1$ has a negative cyclic lift to $\widetilde{\eta} = \sum_{i=1}^{\infty}  \, \eta_i \cdot u^i$. In particular, since $B(\eta_1) = -b(\eta_2)$, if we add $-I(\eta_2)$ to $B( J(s^{n-1}\theta) ) $, then it is $b$-closed. More precisely, we define
\begin{equation}  \label{omega_i_def_fin_cell}
	\begin{split}
		\omega_2 \, &:= \,  B(J(s^{n-1} \theta)) - I(\eta_2)  \, \in \, X^{(2)}(\Pi)_n \\
		\omega_i \, &:= \, - I(\eta_i) \, \in \, X^{(i)}(\Pi)_{n-4+2i} \qquad \text{for } i > 2
	\end{split}
\end{equation}
then one can verify that $\widetilde{\omega} := \sum_{i = 2}^{\infty} \, \omega_i \cdot u^i \, \in \, (F^2 X^{\tot}(\Pi))_{n-4}$ is closed under $d_{\tot} = b + uB$.

We claim that the map $(\omega^2)^{\#}_L : \cSA^{\vee}[n] \ra \cSA$ is still an isomorphism of bimodules in this case. Indeed, consider the following basis of $\cSA$:
\begin{equation*}
	\{ sDf_1, \ldots, sDf_r\} \, \cup \, \{ E_x \}_{x\in \scO} \, \cup \, \{ sDg_1, \ldots, sDg_r\} \, \cup \, \{ sDc_x \}_{x\in \scO} 
\end{equation*}
whose dual basis will be ordered as follows:
\begin{equation*}
(sDg_1)^{\vee}, \ldots, (sDg_r)^{\vee}\} \, \cup \, \{ (sDc_x)^{\vee} \}_{x\in \scO} \cup	\{ (sDf_1)^{\vee}, \ldots, (sDf_r)^{\vee}\} \, \cup \, \{ (E_x)^{\vee} \}_{x\in \scO} 
\end{equation*}

With respect to these basis, we may represent the map $(\omega^2)^{\#}_L$ as a matrix with entries in $\Pi \otimes \Pi$. 
Notice that the formula \eqref{B_J_theta_dual_basis} is still valid in the present deformed case, so that $\omega_2$ is given by \eqref{B_J_theta_dual_basis} minus $I(\eta_2)$. Hence, the matrix that represents $(\omega^2)^{\#}_L$ with respect to the above basis is upper triangular with $\pm 1$ in the diagonals, and hence is an isomorphism.
This completes the proof that $\widetilde{\omega}$ is an $n$-Calabi-Yau structure on $\Pi$.

Now, we give an overview of the proof of the general case (see below). To remove the finitely cellular assumption, several things in the above proof need to be modified.
The first thing is that, since $\cSA$ is no longer finitely cellular (we may and will assume that it is cofibrant), we need to choose a resolution of $\cSA^{\vee}$. In this setting, the strict Casimir element $\theta$ is replaced by a non-strict one (see Lemma \ref{Casimir_element_lemma}). One consequence of this non-strictness is that the equality \eqref{d_J_theta_fin_cell} is replaced by a homotopy equivalence (see \eqref{bh_eta_psi} below). Accordingly, the elements \eqref{omega_i_def_fin_cell} need to be modified by this homotopy (see \eqref{omega_i_def} below). 
Finally, to prove the non-degeneracy of $\omega_2$, we use Lemma \ref{triang_decomp_lemma} below, which may be viewed as an analogue in the general case of the ``upper-triangularity'' of the matrix that represents $(\omega^2)^{\#}_L$.

\bpf[Proof of Theorem \ref{main_thm_abs} (the general case). ]
Assume that $\cA$ is cofibrant%
\footnote{In order to establish the first statement (the existence of Calabi-Yau structure), it suffices to assume that $\cA$ is almost cofibrant. We use cofibrancy only for the second statement (about exactness of the Calabi-Yau structure).}.
Then by Theorem \ref{X_complex_thm_2}, we may represent the deformation parameter $\widetilde{\eta} \in \HN_{n-2}(\cA)$ by a closed element of degree $n-4$ in $F^1X^{\tot}(\cA)$.
In other words, we may write
\begin{equation}  \label{eta_tilde_sum}
\widetilde{\eta} \, = \, \sum_{i=1}^{\infty}  \, \eta_i \cdot u^i
\end{equation}
where $\eta_i \in X^{(i)}(\cA)$ are elements such that 
\begin{equation*}
|\eta_i| = n-4 + 2i \, , \qquad  b(\eta_1) = 0 \, , \qquad b(\eta_i) + B(\eta_{i-1}) = 0 \quad \text{for } i > 1
\end{equation*}

Choose a cofibrant resolution $p : P \xra{\sim} \cSA^{\vee}$. Then an explicit model for the deformed $n$-Calabi-Yau completion is
\begin{equation*}
\Pi \, := \, ( \, T_{\cA}(P[n-1]) \, , \, d_{\eta} \, )
\end{equation*}
where the differential $d_{\eta}$ is induced by the map
\begin{equation}  \label{P_to_A_eta}
P[n-1] \xraq{p} \cSA^{\vee}[n-1] \xraq{(\eta_1)^{\#}} \cA[1]
\end{equation}

Denote by $\iota : \cA \ra \Pi$ the canonical dg functor. This induces a map of $\bN$-graded mixed complexes
\begin{equation}  \label{I_map}
I \, : \, X^{(n)}(\cA) \raq X^{(n)}(\Pi) 
\end{equation} 

We will also be using the following pre-map of complexes (which in general do not commute with the differentials)
\begin{equation}   \label{J_map}
J \, : \, P[n-1] \otimes_{\cAe} \cSA \rintoq T_{\cA}(P[n-1]) \otimes_{\cAe} \cSA \, = \, \Pi \otimes_{\Pie} (\iota^e_!(\cSA)) \xraq{\gamma_{\natural}} \cS(\Pi)_{\natural} \, = \, X^{(1)}(\Pi)
\end{equation} 
where $\gamma = \gamma_{\iota} : \iota^e_!(\cSA) \ra \cS(\Pi)$ is the canonical map.

Notice that \eqref{J_map} does not commute with differential because we used the identification $T_{\cA}(P[n-1]) = \Pi$ in the middle equality. This identification does not commute with the differentials because $T_{\cA}(P[n-1])$ has differential $d_0$ while $\Pi$ has differential $d_{\eta}$. 
In fact, one can compute directly that $d(J)$ is the following composition
\begin{equation}  \label{dJ_map}
d(J) \, : \, 
P[n-1] \otimes_{\cAe} \cSA 
\xraq{\psi}  X^{(1)}(\cA)[1] \xraq{I} X^{(1)}(\Pi)[1]
\end{equation}
where $\psi$ is the map
\begin{equation}  \label{psi_map_1}
	\psi \, : \, P[n-1] \otimes_{\cAe} \cSA 
	\xraq{\eqref{P_to_A_eta} \otimes \id} \cA[1] \otimes_{\cAe} \cSA  \, = \, X^{(1)}(\cA)[1]
\end{equation}

Let $\theta \in Z_0(P \otimes_{\cAe} \cSA)$ be a Casimir element as in Lemma \ref{Casimir_element_lemma}, whose shift will be denoted as $s^{n-1}\theta \in Z_{n-1}( P[n-1] \otimes_{\cAe} \cSA  )$
Applying \eqref{J_map}, we have an element $J(s^{n-1} \theta) \in X^{(1)}(\Pi)_{n-1}$. 
By \eqref{dJ_map}, we see that 
\begin{equation}  \label{b_J_theta}
	b(J(s^{n-1}\theta)) \, = \, I( \psi(s^{n-1}\theta) )
\end{equation}

By Lemma \ref{Casimir_element_lemma}(2), we see that $\eta_1$ is homologous to $\psi(s^{n-1}\theta)$. {\it i.e.,} there exists $h \in X^{(1)}(\cA)_{n-1}$ such that 
\begin{equation}  \label{bh_eta_psi}
	b(h) = \eta_1 - \psi(s^{n-1}\theta)
\end{equation}

We define
\begin{equation}  \label{omega_i_def}
	\begin{split}
	\omega_2 \, &:= \,  B(J(s^{n-1} \theta)) - I(\eta_2) + I(B(h))  \, \in \, X^{(2)}(\Pi)_n \\
    \omega_i \, &:= \, - I(\eta_i) \, \in \, X^{(i)}(\Pi)_{n-4+2i} \qquad \text{for } i > 2
	\end{split}
\end{equation}

We claim that the element $\widetilde{\omega} := \sum_{i = 2}^{\infty} \, \omega_i \cdot u^i \, \in \, (F^2 X^{\tot}(\Pi))_{n-4}$ is an $n$-Calabi-Yau structure on $\Pi$. We first show that it is closed under $d_{\tot}$:
\blm
We have $b(\omega_2) = 0$ and $b(\omega_i) + B(\omega_{i-1}) = 0$ for $i > 2$.
\elm

\bpf
The first statement follows from the computation
\begin{equation}  \label{b_omega_2_zero}
	\begin{split}
		b(\omega_2) \, &= \, b(B(J(s^{n-1}\theta))) - I(b(\eta_2)) + I(b(B(h))) \\
		\, &= \, -B(b(J(s^{n-1}\theta))) + I(B(\eta_1)) - I(B(b(h))) \\
		\, & \stackrel{\eqref{b_J_theta}}{=} \, B( -I(\psi(s^{n-1}\theta)) + I(\eta_1) - I(b(h)) )
		\, \stackrel{\eqref{bh_eta_psi}}{=} \, 0
	\end{split}
\end{equation}
The second statement follows easily from the fact that \eqref{I_map} intertwine with $b$ and $B$.
\epf

By Corollary \ref{TAM_cofib_cor}(2), we see that $\Pi$ is cofibrant, so that by Theorem \ref{X_complex_thm_2} the element $\widetilde{\omega} \in \, Z_{n-4}(F^2 X^{\tot}(\Pi))$ indeed represents a class $[\widetilde{\omega}] \in \HN_n(\Pi)$.
Notice also that $\Pi$ is smooth by Corollary \ref{TAM_cofib_cor}(4).
Thus, to show that $[\widetilde{\omega}]$ is an $n$-Calabi-Yau structure, it remains to show that $\omega_2 \in X^{(2)}(\Pi)$ is non-degenerate. 

By Proposition \ref{Omega_Pi_SES}, there is a short exact sequence in $\Mod(\Pie)$:
\begin{equation} \label{S_Pi_SES_1}
0 \raq \iota^e_! \, \cS(\cA) \xraq{\gamma} \cS(\Pi) \xraq{\rho} \iota^e_! \, P[n] \raq 0
\end{equation}
where $\gamma$ is the usual canonical map, and $\rho$ is the map of bimodule determined by 
\begin{equation}  \label{rho_formula}
	\begin{split}
	\rho(sDf) &= 0 \qquad \text{ for } f \in {}_y\cA_{x} \\
	\rho(E_x) &= 0 \qquad \text{ for } x \in \scO \\
	\rho(sD\xi) &= s \xi \qquad \text{ for } \xi \in {}_y (P[n-1]) _{x}
	\end{split}
\end{equation}

Since the three bimodules in \eqref{S_Pi_SES_1} are cofibrant, they are in particular projective as a graded bimodule. Therefore, taking the duals and shifting by $n$, we still have a short exact sequence
\begin{equation} \label{S_Pi_SES_2}
	0 \raq (\iota^e_! \, P)^{\vee}  \xraq{\rho^{\vee}} \cS(\Pi)^{\vee}[n] \xraq{\gamma^{\vee}}  (\iota^e_! \, \cS(\cA))^{\vee}[n] \raq 0
\end{equation}

Comparing \eqref{S_Pi_SES_1} and \eqref{S_Pi_SES_2}, we see that the first and last terms are quasi-isomorphic. Specifically, denote by $j : P \otimes_{\cAe} \cSA \ra (\iota^e_! \,P) \otimes_{\Pie} (\iota^e_! \,\cSA)$ the canonical map. Then by Lemma \ref{Casimir_element_lemma}(4), $j(\theta) \in Z_0( (\iota^e_! \,P) \otimes_{\Pie} (\iota^e_! \,\cSA) )$ is still a Casimir element, so that it induces quasi-isomorphisms
\begin{equation*}
	\begin{split}
		j(\theta)^{\#}_L \, &: \, (\iota^e_! \, P)^{\vee} \xraq{\sim} \iota^e_! \, \cSA \\
		j(\theta)^{\#}_R \, &: \, (\iota^e_! \, \cSA)^{\vee} \xraq{\sim} \iota^e_! \, P 
	\end{split}
\end{equation*}

The fact that $\omega_2$ is non-degenerate then follows from the following
\blm  \label{triang_decomp_lemma}
Let $\zeta \in X^{(2)}(\cA)_n$ be an element such that $\omega_2 := B(J(s^{n-1} \theta)) + I(\zeta) \in X^{(2)}(\Pi)_n$ is $b$-closed, then the following diagram in $\Mod(\Pie)$ commutes:
\begin{equation}  \label{triang_decomp_comm_diag}
	\begin{tikzcd}
		0 \ar[r] & (\iota^e_! \, P)^{\vee}  \ar[r,"\rho^{\vee}"] \ar[d, "j(\theta)^{\#}_L"] & \cS(\Pi)^{\vee}[n] \ar[r, "\gamma^{\vee}"] \ar[d, "(\omega_2)^{\#}_L"] & (\iota^e_! \, \cS(\cA))^{\vee}[n] \ar[r] \ar[d, "j(\theta)^{\#}_R"] & 0 \\
		0 \ar[r]  & \iota^e_! \, \cS(\cA) \ar[r, "\gamma"] & \cS(\Pi) \ar[r, "\rho"] & \iota^e_! \, P[n] \ar[r] & 0
	\end{tikzcd}
\end{equation}
\elm

\bpf
The pre-map $I(\zeta)^{\#}_L : \cS(\Pi)^{\vee}[n] \ra \cS(\Pi)$ induced by $I(\zeta)$ is given by the composition
\begin{equation*}
	\cS(\Pi)^{\vee}[n] \xraq{\gamma^{\vee}} (\iota^e_! \, \cS(\cA))^{\vee}[n] 
	\xraq{j(\zeta)^{\#}_L } \iota^e_! \, \cS(\cA) \xraq{\gamma}  \cS(\Pi)
\end{equation*} 
Therefore, both its pre-composition with $(\iota^e_! \, P)^{\vee} \xra{\rho^{\vee}} \cS(\Pi)^{\vee}[n]$ and its post-composition with $\cS(\Pi) \xra{\rho} \iota^e_! \, P[n]$ are zero. 
Thus, to verify the commutativity of \eqref{triang_decomp_comm_diag}, we may assume that $\zeta = 0$, except that in this case, the middle vertical arrow $(\omega_2)^{\#}_L$ in \eqref{triang_decomp_comm_diag} is only a pre-map. But this is not problematic because we can still discuss the commutativity of a diagram of pre-maps.

We now show that, for any element $\theta \in (P \otimes_{\cAe} \cSA)_0$ (not necessarily a Casimir element, not even necessarily closed), the diagram of pre-maps \eqref{triang_decomp_comm_diag} commutes, where we take $\omega_2 := B(J(s^{n-1} \theta))$.
Since we assert it for any element $\theta$, it suffices to verify the statement for the elements of the two types: 
\begin{enumerate}
	\item $\theta = \xi \otimes sDf$ for $f \in \cA(x,y)_{m-1}$ and $\xi \in P(y,x)_{-m}$.
	\item $\theta = \xi \otimes E_x$ for $x \in \scO$ and $\xi \in P(x,x)_0$.
\end{enumerate}

For both types, we may compute $\omega_2$ by Lemma \ref{B_on_X1}. For type (1), we have $\omega_2 = sD(s^{n-1}\xi) \otimes sDf + (-1)^{m(n+m)} sDf \otimes sD(s^{n-1}\xi)$. For type (2), we have $\omega_2 = sD(s^{n-1}\xi) \otimes E_x + E_x \otimes sD(s^{n-1}\xi)$.
The commutativity of \eqref{triang_decomp_comm_diag} is then a straightforward verification for both cases.
\epf

This completes the proof that $[\widetilde{\omega}]$ is an $n$-Calabi-Yau structure. Thus, it suffices to verify that $[\widetilde{\omega}]$ is exact if $[\widetilde{\eta}]$ is exact. 
Notice that the image of $B : \HC_{m-1} \ra \HN_m$ coincides with the image of $\overline{B} : \overline{\HC}_{m-1} \ra \HN_m$. Thus, to work with exactness of $[\widetilde{\eta}]$ and $[\widetilde{\omega}]$, we may use Theorem \ref{X_complex_thm_3}. 
More precisely, we apply Theorem \ref{X_complex_thm_3} for $r = 1$, so that we may assume that the element $\widetilde{\eta}$ in \eqref{eta_tilde_sum} is given by $u \cdot \overline{B}(\zeta_0)$ for some $\zeta_0 \in Z_{n-3}(\overline{\cA_{\natural}})$. In other words, we have $\eta_i = 0$ for $i > 1$, and $\eta_1 = \overline{B}(\zeta_0)$.
Then we will apply Theorem \ref{X_complex_thm_3} for $r = 2$ to $\Pi$, so that it suffices to show that $\widetilde{\omega} = u^2\cdot B(\lambda_1)$ for some closed element $\lambda_0 + u \cdot \lambda_1 \in \Tot \, X^{<2}_{{\rm sr}}(\Pi)$. 
In other words, we will construct
\begin{equation*}
	\lambda_0 \in (\overline{\Pi_{\natural}})_{n-3}  \qquad \text{and} \qquad \lambda_1 \in X^{(1)}(\Pi)_{n-1}
\end{equation*}
such that
\begin{equation}  \label{omega_lambda_relations}
	b(\lambda_0) = 0 \, , \qquad  b(\lambda_1) + \overline{B}(\lambda_0) = 0 \, , \qquad \omega_2 = B(\lambda_1) \, , \qquad \omega_i = 0 \quad \text{for } i > 2
\end{equation}

Denote by $I : \overline{\cA_{\natural}} \ra \overline{\Pi_{\natural}}$ the canonical map induced by $\iota : \cA \ra \Pi$. Thus, together with \eqref{I_map}, it gives a map of $\bN$-graded mixed complexes $I : X^{\bullet}_{{\rm sr}}(\cA) \ra X^{\bullet}_{{\rm sr}}(\Pi)$. Define
\begin{equation*}
	\lambda_0 := -I(\zeta_0) \qquad \text{and} \qquad \lambda_1 := J(s^{n-1} \theta) + I(h)
\end{equation*}
Then the last two equations of \eqref{omega_lambda_relations} follows directly from the definition \eqref{omega_i_def} of $\omega_i$. The equation $b(\lambda_0) = 0$ follows from $b(\zeta_0) = 0$. Hence, it suffices to verify
\begin{equation*}
	\begin{split}
		b(\lambda_1) + \overline{B}(\lambda_0) \, &= \, b(J(s^{n-1} \theta)) + I(b(h)) - \overline{B}(I(\zeta_0)) \\
		\, & \stackrel{\eqref{b_J_theta}}{=} \, I(\psi(s^{n-1}\theta) + b(h) - \overline{B}(\zeta_0))  \, \stackrel{\eqref{bh_eta_psi}}{=} \, 0
	\end{split}
\end{equation*} 
\epf

\brm  \label{CY_str_can_remark}
Our definition \eqref{omega_i_def} of the $n$-Calabi-Yau structure on $\Pi$ depends on a choice of the homotopy $h \in X^{(1)}(\cA)_{n-1}$. Any other such choice differs by an element in $Z_{n-1}( X^{(1)}(\cA) )$. Thus, our construction of the $n$-Calabi-Yau structure $[\widetilde{\omega}] \in \HN_{n}(\Pi)$ is ambiguous up to addition of an element in the image of the map $\HH_{n-1}(\cA)  \xra{B} \HN_n(\cA) \xra{I} \HN_{n}(\Pi)$. Indeed, by Lemma \ref{triang_decomp_lemma}, $\widetilde{\omega}$ remains non-degenerate after adding any such element. Notice also that, the element $h$ is constructed in the proof of Lemma \ref{Casimir_element_lemma}(2), so that in particular, it depends on the given homotopy $\varphi' \simeq \id$ in Lemma \ref{Casimir_element_lemma}(1). 
Thus, this (non-)canonicality of $[\widetilde{\omega}]$ stems from the (non-)canonicality of the ``Casimir data'' $(\theta, \varphi' \simeq \id)$ of Lemma \ref{Casimir_element_lemma}(1). It is not clear to us how (and to what extent) one can formulate a sharper canonicality of this ``Casimir data'', and how these ``Casimir data'' relate to each other under derived equivalences $\cA' \xra{\sim} \cA$ and under change of choices $p : \cP \xra{\sim} \cSA^{\vee}$. 
\erm

\brm
Our proof also shows that, if we are given $[\eta] \in \HH_{n-2}(\cA)$ that is weakly closed, in the sense that its image under $B : \HH_{n-2}(\cA) \ra \HH_{n-1}(\cA)$ is zero, then the ``deformed Calabi-Yau completion'' $\Pi_n(\cA,\eta)$ is weakly $n$-Calabi-Yau, meaning that there is an isomorphism $\Pi^![n] \cong \Pi$ in $\cD(\Pie)$. Namely, choose $\eta_1 \in Z_{n-2}(X^{(1)}(\cA))$ that represents $[\eta]$, then the weakly closed condition guarantees that $B(\eta_1) + b(\eta_2) = 0$ for some $\eta_2 \in X^{(1)}(\cA)_{n}$. This data is enough to define $\omega_2$ in \eqref{omega_i_def}, and to verify that $b(\omega_2) = 0$ in \eqref{b_omega_2_zero}. The proof that $\omega_2$ is non-degenerate goes though unchanged.
\erm

\section{Main theorem: the relative case}  \label{sec_main_thm_rel}

\bthm  \label{main_thm_rel}
For any dg functor $F: \cA \ra \cB$ between smooth $k$-flat small dg categories, and any relative negative cyclic class $\widetilde{\eta} \in \HN_{n-2}(\cB,\cA)$, the deformed relative $n$-Calabi-Yau completion $\Pi_{n-1}(\cA, \eta_{\cA}) \raq \Pi_n(\cB,\cA, \eta)$ has a relative $n$-Calabi-Yau structure.

Moreover, if the deformation parameter $\widetilde{\eta}$ is exact in the sense that it is in the image of $B : \HC_{n-3}(\cB,\cA) \ra \HN_{n-2}(\cB,\cA)$, then this relative $n$-Calabi-Yau structure on $\Pi_n(\cB,\cA,\eta)$ is exact.
\ethm

\bpf
We first discuss Casimir elements in the relative setting. The natural setting for this discussion starts with $Q_{\cA} \in \Mod(\cAe)$ and $Q_{\cB} \in \Mod(\cBe)$, both cofibrant, and a map $\gamma : F^e_!(Q_{\cA}) \ra Q_{\cB}$ in $\Mod(\cBe)$. However, we will only need the case $Q_{\cA} = \cSA$ and $Q_{\cB} = \cSB$ and $\gamma = \gamma_F$ being the natural map (implicitly this means we will require $\cA$ and $\cB$ to be almost cofibrant), so we will focus on this case.

Assume that $\cA$ and $\cB$ are cofibrant. 
Choose a cofibrant replacement $p_{\cA} : P_{\cA} \xra{\sim} \cSA^{\vee}$ in $\Mod(\cAe)$ and a cofibrant replacement $p_{\cB} : P_{\cB} \xra{\sim} \cSB^{\vee}$ in $\Mod(\cBe)$.
Recall from \eqref{F_shriek_dual_2} that the map 
\begin{equation*}
	q_{\cA} \, : \, F^e_!(P_{\cA}) \xraq{F^e_!(p_{\cA})} F^e_!(\cSA^{\vee}) \raq F^e_!(\cSA)^{\vee}
\end{equation*}
is a quasi-isomorphism. Choose a map $\nu : P_{\cB} \ra F^e_!(P_{\cA})$ that represents the map $\gamma_F^{\vee} : \cSB^{\vee} \ra F^e_!(\cSA)^{\vee}$ in the derived category $\cD(\cBe)$. In other words, the following diagram commutes up to homotopy:
\begin{equation}  \label{P_BA_cS_BA_homotopy_comm}
	\begin{tikzcd}
		P_{\cB} \ar[r, "\nu"] \ar[d, "p_{\cB}"', "\sim"] & F^e_!(P_{\cA}) \ar[d, "q_{\cA}", "\sim"'] \\
		\cSB^{\vee} \ar[r, "\gamma_F^{\vee}"] & F^e_!(\cSA)^{\vee}
	\end{tikzcd}
\end{equation}

Choose such a homotopy $\chi \in \Homcom_{\cB}(P_{\cB},F^e_!(\cSA)^{\vee})_1$ between $\gamma_F^{\vee} \circ p_{\cB}$ and $q_{\cA} \circ \nu$. 
Denote by
\begin{equation}  \label{PBA_SBA_def}
	\begin{split}
		P_{\cB/\cA} \, &:= \, \cone[P_{\cB} \xra{\nu} F^e_!(P_{\cA})][-1] \\
		\cS(\cB,\cA) \, &:= \, \cone[F^e_!(\cSA) \xra{\gamma_F} \cSB]
	\end{split}
\end{equation}
Then the maps ($p_{\cB},q_{\cA},\chi$) defines a map in $\Mod(\cBe)$:
\begin{equation}  \label{p_BA_def}
	p_{\cB/\cA} \, : \, P_{\cB/\cA} \raq \cS(\cB,\cA)^{\vee}
\end{equation}

Clearly, the following diagram in $\Mod(\cBe)$ commutes:
\begin{equation}  \label{PBA_SBA}
	\begin{tikzcd}
		 F^e_!(P_{\cA})[-1] \ar[r] \ar[d, "q_{\cA}"]  &  P_{\cB/\cA} \ar[d, "p_{\cB/\cA}"] \ar[r] & P_{\cB} \ar[d, "p_{\cB}"]  \\
		 F^e_!(\cSA)^{\vee}[-1] \ar[r] &  \cS(\cB,\cA)^{\vee} \ar[r] & \cSB^{\vee} 
	\end{tikzcd}
\end{equation}
where the horizontal maps are the canonical maps relating the cone to its constituents. Since the connecting square \eqref{P_BA_cS_BA_homotopy_comm} is commutative up to homotopy, we see that \eqref{PBA_SBA} is a map of exact triangles in $\cD(\cBe)$. Since $q_{\cA}$ and $p_{\cB}$ are quasi-isomorphisms, we see therefore that $p_{\cB/\cA}$ is also a quasi-isomorphism.

The canonical map $\widetilde{\mu} : F^e_!(P_{\cA})[-1] \ra  P_{\cB/\cA}$ in $\Mod(\cBe)$ corresponds by adjunction $F^e_! \dashv F^{e*}$  to a map
\begin{equation}
	\mu \, : \, P_{\cA}[-1] \raq P_{\cB/\cA} \qquad \text{ in } \Mod(\cAe)
\end{equation}

Now suppose we are given $M_{\cA} \in \Mod(\cAe)$ and $M_{\cB} \in \Mod(\cBe)$, together with a map $\mu_M : M_{\cA} \ra M_{\cB}$ in $\Mod(\cAe)$. Then we consider the complex
\begin{equation*}
	\Xi(M_{\cB}, M_{\cA}) \, := \, \cone\, [ \, M_{\cA} \otimes_{\cAe} \cSA \xraq{\mu_M \otimes \gamma_F} M_{\cB} \otimes_{\cBe} \cSB \, ]
\end{equation*}

By the construction in the proof of Proposition \ref{rel_Hoch_induce_map_of_triag_prop}, we see in particular that any homogeneous element $\xi \in \Xi(M_{\cB},M_{\cA})_m$ determines a commutative diagram of pre-maps (of degree $0$) in $\underline{\Mod}_0(\cAe)$:
\begin{equation}   \label{SBA_MBA_xi}
	\begin{tikzcd}
		F^e_!(\cSA)^{\vee}[m-1] \ar[r]  \ar[d, "(\xi'')^{\#}_R"] &  \cS(\cB,\cA)^{\vee}[m]  \ar[d, "(\xi')^{\#}_R"] \\
		F^e_!(M_{\cA}) \ar[r] & M_{\cB}
	\end{tikzcd}
\end{equation}

Denote by $\xi_{\cA} \in (M_{\cA} \otimes_{\cAe} \cSA)_{m-1}$ the image of $\xi$ under the canonical map $\Xi(M_{\cB},M_{\cA}) \ra M_{\cA} \otimes_{\cAe} \cSA[1]$. Then $\xi_{\cA}$ induces a pre-map in $\underline{\Mod}_0(\cAe)$:
\begin{equation}  \label{PA_MA_xi}
	P_{\cA}[m-1] \xraq[\sim]{p_{\cA}} \cSA^{\vee}[m-1] \xraq{(\xi_{\cA})^{\#}_R} M_{\cA}
\end{equation}

On the other hand, there is a pre-map in $\underline{\Mod}_0(\cBe)$:
\begin{equation}  \label{PBA_MBA_xi}
	P_{\cB/\cA}[m] \xraq[\sim]{p_{\cB/\cA}} \cS(\cB,\cA)^{\vee}[m] \xraq{(\xi')^{\#}_R} M_{\cB}
\end{equation}

Combine \eqref{SBA_MBA_xi} with the left hand side square of \eqref{PBA_SBA}, and noting that the composition $(\xi'')^{\#}_R \circ q_{\cA}$ is precisely the image  of $\eqref{PA_MA_xi}$ under $F^e_!$, we see that the following diagram in $\underline{\Mod}_0(\cBe)$ commutes:
\begin{equation}  \label{xi_PA_PAB_comm_1}
	\begin{tikzcd}
		F^e_!(P_{\cA})[m-1] \ar[r] \ar[d, "F^e_!(\eqref{PA_MA_xi})"'] &  P_{\cB/\cA}[m] \ar[d, "\eqref{PBA_MBA_xi}"] \\
		F^e_!(M_{\cA}) \ar[r] &  M_{\cB}
	\end{tikzcd}
\end{equation}
By adjunction, this can be rewritten as a commutative diagram in $\underline{\Mod}_0(\cAe)$:
\begin{equation}  \label{xi_PA_PAB_comm_2}
	\begin{tikzcd}
	P_{\cA}[m-1] \ar[r, "\mu"] \ar[d, "\eqref{PA_MA_xi}"'] &  P_{\cB/\cA}[m] \ar[d, "\eqref{PBA_MBA_xi}"] \\
M_{\cA} \ar[r, "\mu_M"] &  M_{\cB}
	\end{tikzcd}	
\end{equation}

Denote by $\Homcom_{\cBe/\cAe}(\mu,\mu_M)$ the complex whose degree $m$ elements are
\begin{equation*}
	\Homcom_{\cBe/\cAe}(\mu,\mu_M)_m \, := \, \left \{ \, (f,g) \, \middle | \, \parbox{14em}{$f \in \Homcom_{\cA}(P_{\cA}[-1],M_{\cA})_m$ \\ 
	$g \in  \Homcom_{\cB}(P_{\cB/\cA},M_{\cB})_m$ \\
such that the following diagram of pre-maps commutes:} \, \, \,  \begin{tikzcd}
P_{\cA}[m-1] \ar[r, "\mu"] \ar[d, "f"'] &  P_{\cB/\cA}[m] \ar[d, "g"] \\
M_{\cA} \ar[r, "\mu_M"] &  M_{\cB}
\end{tikzcd}  \right \}
\end{equation*}
and whose differential is defined by $d(f,g) := (d(f),d(g))$.

Then our above discussion defines a map of chain complexes
\begin{equation} \label{Theta_to_Hom_BA}
	\Xi(M_{\cB},M_{\cA}) \raq \Homcom_{\cBe/\cAe}(\mu,\mu_M)
\end{equation}

\blm  \label{Theta_to_Hom_BA_qism}
The map \eqref{Theta_to_Hom_BA} is a quasi-isomorphism.
\elm

\bpf
Since $P_{\cB/\cA}$ is defined as a cone \eqref{PBA_SBA_def}, a pre-map $g \in \Homcom_{\cBe}(P_{\cB/\cA},M_{\cB})_m$ is given by $(g_1,g_2)$ where $g_1 \in \Homcom_{\cBe}(F^e_!(P_{\cA}[-1]),M_{\cB})_m$ and $g_2 \in \Homcom_{\cBe}(P_{\cB},M_{\cB})_m$. 
To specify an element $(f,g) \in \Homcom_{\cBe/\cAe}(\mu,\mu_M)_m$, once $f$ is specified, the component $g_1$ of $g$ is uniquely determined by the commutativity requirement $g \circ \mu = \mu_M \circ f$. Thus, the pair $(f,g)$ is uniquely determined by $(f,g_2)$, with no constraints.
In other words, negelecting differentials, we have an isomorphism
\begin{equation*}
	\Homcom_{\cBe/\cAe}(\mu,\mu_M) \, \cong \, \Homcom_{\cAe}(P_{\cA}[-1], M_{\cA}) \, \oplus \, \Homcom_{\cBe}(P_{\cB},M_{\cB})
\end{equation*}

Putting in the differentials, we see that $\Homcom_{\cBe/\cAe}(\mu,\mu_M)$ is the cone of the following map:
\begin{equation}  \label{mu_nu_star}
	\Homcom_{\cAe}(P_{\cA}, M_{\cA}) \xraq{(\mu_M)_*} \Homcom_{\cAe}(P_{\cA}, M_{\cB}) \, = \, \Homcom_{\cBe}(F^e_!(P_{\cA}), M_{\cB}) \xraq{\nu^*} \Homcom_{\cBe}(P_{\cB}, M_{\cB})
\end{equation}

To specify the map \eqref{Theta_to_Hom_BA} under this identification, consider the following diagram:
\begin{equation}  \label{MSA_Hom_PMA_square}
	\begin{tikzcd}
		M_{\cA} \otimes_{\cAe} \cSA \ar[r, "\mu_M \otimes \gamma_F"] \ar[d] \ar[rd, dotted] & M_{\cB} \otimes_{\cBe} \cSB \ar[d] \\
		\Homcom_{\cAe}(P_{\cA}, M_{\cA}) \ar[r, "\eqref{mu_nu_star}"] & \Homcom_{\cBe}(P_{\cB}, M_{\cB})
	\end{tikzcd}
\end{equation}
where the vertical ones are the canonical ones (induced by $p_{\cA}$ and $p_{\cB}$), and the dotted arrow is a homotopy between the two compositions, and is defined by
\begin{equation*}
	M_{\cA} \otimes_{\cAe} \cSA \xraq{(\mu_M)_*} \Homcom_{\cB}(F_e^!(\cSA)^{\vee}, M_{\cB}) \xraq{\chi^*} \Homcom_{\cB}(P_{\cB}, M_{\cB})[-1]
\end{equation*}

The two vertical maps and the dotted homotopy in \eqref{MSA_Hom_PMA_square} together defines a map on the cones of the horizontal maps in \eqref{MSA_Hom_PMA_square}. Under the above identification of $\Homcom_{\cBe/\cAe}(\mu,\mu_M)$ with the cone of \eqref{mu_nu_star}, one can verify that this is precisely the map \eqref{Theta_to_Hom_BA} (recall that the pre-map $\chi$ was used in the definition of \eqref{p_BA_def}, and hence it implicitly appears in the map \eqref{Theta_to_Hom_BA}).
Therefore, \eqref{Theta_to_Hom_BA} fits into a map of exact triangles together with the two vertical maps of \eqref{MSA_Hom_PMA_square}. Since both of the two vertical maps of \eqref{MSA_Hom_PMA_square} are quasi-isomorphisms, so is \eqref{Theta_to_Hom_BA}.
\epf

The quasi-isomorphism \eqref{Theta_to_Hom_BA} allows us to establish the basic properties of ``relative Casimir elements''. For this, we take $M_{\cA} = P_{\cA}[-1]$, $M_{\cB} = P_{\cB/\cA}$, and $\mu_M = \mu$. In this case, we denote by
\begin{equation*}
	\Theta_{\cB/\cA} \, := \,	\Xi(P_{\cB/\cA},P_{\cA}[-1])
\end{equation*}

\bcor  \label{rel_Casimir_cor}
\begin{enumerate}
	\item There exists an element $\theta_{\cB/\cA} \in Z_0(\Theta_{\cB/\cA})$ whose image under the map $\Theta_{\cB/\cA} \ra \Homcom_{\cBe/\cAe}(\mu,\mu)$ defined in \eqref{Theta_to_Hom_BA} is homologous to the identity. We call such an element a \emph{relative Casimir element}.
	\item Let $\xi \in Z_m( \Xi(M_{\cB},M_{\cA}) )$ be any closed element, which therefore induces the commutative diagram \eqref{xi_PA_PAB_comm_2} in $\Mod(\cAe)$, and hence a commutative diagram of chain complexes
	\begin{equation*}
		\begin{tikzcd}
			P_{\cA}[m-1] \otimes_{\cAe} \cSA \ar[r, "\mu \otimes \gamma_F"] \ar[d, "\eqref{PA_MA_xi} \otimes \id"'] &  P_{\cB/\cA}[m]  \otimes_{\cBe} \cSB \ar[d, "\eqref{PBA_MBA_xi} \otimes \id"] \\
			M_{\cA}  \otimes_{\cAe} \cSA \ar[r, "\mu_M \otimes \gamma_F"] &  M_{\cB} \otimes_{\cBe} \cSB 
		\end{tikzcd}
	\end{equation*}
Taking the cone of the horizontal maps, we obtain an induced map $\psi : \Theta(\cB,\cA)[m] \ra \Xi(M_{\cB},M_{\cA})$. Then for any relative Casimir element $\theta_{\cB/\cA} \in Z_0(\Theta_{\cB/\cA})$, its image $\psi(s^m \theta_{\cB/\cA})$ is homologous to $\xi$.
\end{enumerate}
\ecor

\bpf
(1) clearly follows from Lemma \ref{Theta_to_Hom_BA_qism}. For (2), it suffices to show that the two elements $\psi(\theta_{\cB/\cA})$ and $\xi$ have homologous images under \eqref{Theta_to_Hom_BA}. Clearly, the map \eqref{Theta_to_Hom_BA} is natural with respect to maps of $(M_{\cA},M_{\cB},\mu_M)$. Thus, the image of $\psi(\theta_{\cB/\cA})$ under \eqref{Theta_to_Hom_BA} is given by the composition
\begin{equation*}
	\begin{tikzcd}
		P_{\cA}[m-1] \ar[r] \ar[d, "\mu"] & P_{\cA}[m-1] \ar[r] \ar[d, "\mu"] & M_{\cA} \ar[d, "\mu_M"] \\
		P_{\cB/\cA}[m] \ar[r] & P_{\cB/\cA}[m] \ar[r] & M_{\cB}
	\end{tikzcd}
\end{equation*}
where the left hand side square is the map in $\Homcom_{\cBe/\cAe}(\mu,\mu)$ induced by $\theta_{\cB/\cA}$, while the right hand side square is the map in $\Homcom_{\cBe/\cAe}(\mu,\mu_M)$ induced by $\xi$. Therefore the composition is homologous to the map in $\Homcom_{\cBe/\cAe}(\mu,\mu_M)$ induced by $\xi$.
\epf

We now proceed to give an explicit model for the deformed relative $n$-Calabi-Yau completion. In fact, the above discussion will allow us to specify the differential induced by the deformation parameter.

Define $X^{(n)}(\cB, \cA) := \cone[X^{(n)}(\cA) \ra X^{(n)}(\cB)]$, which again forms an $\bN$-graded mixed complex under the induced maps $b$ and $B$. Represent the deformation parameter $\widetilde{\eta} \in \HN_{n-2}(\cB,\cA)$ by a closed element in $F^1 X^{\tot}(\cB , \cA)$. In other words, we write
\begin{equation*}  
	\widetilde{\eta} \, = \, \sum_{i=1}^{\infty}  \, \eta_i \cdot u^i
\end{equation*}
where $\eta_i \in X^{(i)}(\cB,\cA)$ are elements such that 
\begin{equation*}
	|\eta_i| = n-4 + 2i \, , \qquad  b(\eta_1) = 0 \, , \qquad b(\eta_i) + B(\eta_{i-1}) = 0 \quad \text{for } i > 1
\end{equation*}

Notice that $X^{(1)}(\cB,\cA)$ is precisely the complex $X^{(1)}(\cB,\cA) = \Xi(M_{\cB},M_{\cA})$ for the data $M_{\cA} = \cA$, $M_{\cB} = \cB$ and $\mu_M = F$. In particular, applying the map \eqref{Theta_to_Hom_BA} to the lowest order element $\eta_1 \in Z_{n-2}(X^{(1)}(\cB,\cA))$ of $\widetilde{\eta}$, we obtain a map $(\delta_{\eta_{\cA}},\delta_{\eta}) \in Z_{n-2}( \Homcom_{\cBe/\cAe}( \mu , F ) )$. We write this as a commutative diagram
\begin{equation}  \label{eta_PA_PAB_comm}
	\begin{tikzcd}
		P_{\cA}[n-2] \ar[r, "\mu"] \ar[d, "\delta_{\eta_{\cA}}"'] &  P_{\cB/\cA}[n-1] \ar[d, "\delta_{\eta}"] \\
		\cA[1] \ar[r, "F"] &  \cB[1]
	\end{tikzcd}
\end{equation}

As a consequence, we have a commutative diagram of chain complexes
\begin{equation}  \label{psi_comm_diag}
	\begin{tikzcd}
		P_{\cA} \otimes_{\cAe} \cSA [n-2] \ar[r, "\mu \otimes \gamma_F"] \ar[d, "\delta_{\eta_{\cA}} \otimes \id"'] &  P_{\cB/\cA} \otimes_{\cBe} \cSB[n-1] \ar[d, "\delta_{\eta} \otimes \id"] \\
		\cA \otimes_{\cAe} \cSA [1] \ar[r, "F \otimes \gamma_F"] &  \cB \otimes_{\cBe} \cSB[1]
	\end{tikzcd}
\end{equation}
Taking the cone of the horizontal maps, we have an induced map
\begin{equation}  \label{psi_Theta_X}
	\psi \, : \, \Theta_{\cB/\cA}[n-1] \raq X^{(1)}(\cB,\cA)[1]
\end{equation}
By Corollary \ref{rel_Casimir_cor}, we see that $\psi(s^{n-1}\theta_{\cB/\cA})$ is homotopic to $\eta_1$. {\it i.e.}, there exists $h \in X^{(1)}(\cB/\cA)_{n-1}$ such that
\begin{equation}  \label{bh_eta_psi_rel}
	b(h) = \eta_1 - \psi(s^{n-1}\theta_{\cB/\cA})
\end{equation}

Define
\begin{equation*}
	\begin{split}
		\Pi_{\cA} \, &:= \, T_{\cA}( P_{\cA}[n-2] ) \, , \qquad d = d_{\eta_{\cA}} \text{ induced by } \delta_{\eta_{\cA}} \text{ in } \eqref{eta_PA_PAB_comm}  \\
		\Pi_{\cB/\cA} \, &:= \, T_{\cB}( P_{\cB/\cA}[n-1] ) \, , \qquad d = d_{\eta} \text{ induced by } \delta_{\eta} \text{ in } \eqref{eta_PA_PAB_comm}
	\end{split}
\end{equation*}
Then by the commutative diagram \eqref{eta_PA_PAB_comm}, one has a dg functor
\begin{equation*}
	\widetilde{F} \ : \, \Pi_{\cA} \raq \Pi_{\cB/\cA}
\end{equation*}
given by $\widetilde{F}|_{\cA} = F$ and $\widetilde{F}|_{P_{\cA}[n-2]} = \mu$.
The dg functor $\widetilde{F}$ will be taken as our explicit model for the deformed relative $n$-Calabi-Yau completion.

Recall the pre-map \eqref{J_map}. There are obvious analogues for $\Pi_{\cA}$ and $\Pi_{\cB/\cA}$, which define the vertical pre-maps of the following diagram:
\begin{equation*}
	\begin{tikzcd}
		P_{\cA}[n-2] \otimes_{\cAe} \cSA  \ar[r, "\mu \otimes \gamma_F"] \ar[d, "J"]
		&  P_{\cB/\cA}[n-1] \otimes_{\cBe} \cSB \ar[d, "J"]  \\
		X^{(1)}(\Pi_{\cA}) \ar[r]
		& X^{(1)}(\Pi_{\cB/\cA}) 
	\end{tikzcd}
\end{equation*}

This diagram of pre-maps is clearly commutative. Taking the cone of the horizontal maps, we obtain a pre-map
\begin{equation}  \label{J_map_AB}
	J \, : \, \Theta_{\cB/\cA}[n-1] \raq X^{(1)}(\Pi_{\cB/\cA},\Pi_{\cA})
\end{equation}

By a computation parallel to the absolute case, the differential $d(J)$ can be shown to be given by the composition
\begin{equation*}
	d(J) \, : \, \Theta_{\cB/\cA}[n-1] \xraq{\psi} X^{(1)}(\cB,\cA)[1] \xraq{I} X^{(1)}(\Pi_{\cB/\cA} , \Pi_{\cA})[1]
\end{equation*}
where $\psi$ is the map defined in \eqref{psi_Theta_X}, and $I$ is the canonical map of $\bN$-graded mixed complexes
\begin{equation*}
	I \, : \, X^{(m)}(\cB,\cA) \raq X^{(m)}(\Pi_{\cB/\cA} , \Pi_{\cA})
\end{equation*}
induced by the canonical inclusions $\iota_{\cA} : \cA \ra \Pi_{\cA}$ and $\iota_{\cB} : \cB \ra \Pi_{\cB/\cA}$.

In particular, we have
\begin{equation}  \label{b_J_theta_rel}
	b(J(s^{n-1}\theta_{\cB/\cA})) \, = \, I(\psi(s^{n-1}\theta_{\cB/\cA}))
\end{equation}

We define
\begin{equation}  \label{omega_i_def_rel}
	\begin{split}
		\omega_2 \, &:= \,  B(J(s^{n-1} \theta_{\cB/\cA})) - I(\eta_2) + I(B(h))  \, \in \, X^{(2)}(\Pi_{\cB/\cA},\Pi_{\cA})_n \\
		\omega_i \, &:= \, - I(\eta_i) \, \in \, X^{(i)}(\Pi_{\cB/\cA},\Pi_{\cA})_{n-4+2i} \qquad \text{for } i > 2
	\end{split}
\end{equation}

We claim that the element $\widetilde{\omega} := \sum_{i = 2}^{\infty} \, \omega_i \cdot u^i \, \in \, (F^2 X^{\tot}(\Pi_{\cB/\cA},\Pi_{\cA}))_{n-4}$ is a relative $n$-Calabi-Yau structure on $\widetilde{F} : \Pi_{\cA} \ra \Pi_{\cB/\cA}$. We first show that it is closed under $d_{\tot}$:
\blm
We have $b(\omega_2) = 0$ and $b(\omega_i) + B(\omega_{i-1}) = 0$ for $i > 2$.
\elm

\bpf
This follows from the exact same computation as in \eqref{b_omega_2_zero}, replacing \eqref{b_J_theta} and \eqref{bh_eta_psi} by \eqref{b_J_theta_rel} and \eqref{bh_eta_psi_rel} respectively.
\epf

We now prove that $\omega_2$ is non-degenerate. First notice that $(\omega_2)_{\Pi_{\cA}} \in X^{(2)}(\Pi_{\cA})_{n-1}$ is precisely the one constructed in the proof of Theorem \ref{main_thm_abs}, and hence is non-degenerate. By the discussion in the paragraph that follows Definition \ref{rel_CY_def}, it suffices to verify that either $(\omega_2')^{\#}_L$ or $(\omega_2''')^{\#}_L$ in the following induced diagram is a quasi-isomorphism:
\begin{equation*} 
	\begin{tikzcd}
		\cS(\Pi_{\cB/\cA})^{\vee}[n-1] \ar[r] \ar[d, "(\omega_2')^{\#}_L"]  & \widetilde{F}^e_!(\cS(\Pi_{\cA}))^{\vee}[n-1] \ar[r]  \ar[d, "(\omega_2'')^{\#}_L"] &  \cS(\Pi_{\cB/\cA}, \Pi_{\cA})^{\vee}[n]  \ar[d, "(\omega_2''')^{\#}_L"]\\
		\cS(\Pi_{\cB/\cA}, \Pi_{\cA})[-1] \ar[r] & F^e_!(\cS(\Pi_{\cA})) \ar[r] & \cS(\Pi_{\cB/\cA})
	\end{tikzcd}
\end{equation*}

By Proposition \ref{Omega_Pi_SES}, there are short exact sequences
\begin{equation}  \label{S_Pi_A_SES}
	0 \raq (\iota_{\cA})^e_! \, \cSA \xraq{\gamma} \cS(\Pi_{\cA}) \xraq{\rho} (\iota_{\cA})^e_! \, P_{\cA}[n-1] \raq 0
\end{equation}
\begin{equation}  \label{S_Pi_BA_SES}
	0 \raq (\iota_{\cB})^e_! \, \cSB \xraq{\gamma} \cS(\Pi_{\cB/\cA}) \xraq{\rho} (\iota_{\cB})^e_! \, P_{\cB/\cA}[n] \raq 0
\end{equation}

There is an obvious map from $\widetilde{F}^e_!(\eqref{S_Pi_A_SES})$ to \eqref{S_Pi_BA_SES}:
\begin{equation}  \label{S_Pi_SES_A_to_BA}
	\begin{tikzcd}
		0 \ar[r] & (\iota_{\cB})^e_! F^e_! \, \cSA \ar[r, "\gamma"] \ar[d, "\gamma"] & \widetilde{F}^e_! \cS(\Pi_{\cA}) \ar[r, "\rho"] \ar[d, "\gamma"] & (\iota_{\cB})^e_! F^e_! \, P_{\cA}[n-1] \ar[r] \ar[d, "(\iota_{\cB})^e_!(\widetilde{\mu})"] & 0 \\
		0 \ar[r] & (\iota_{\cB})^e_! \, \cSB \ar[r, "\gamma"] & \cS(\Pi_{\cB/\cA}) \ar[r, "\rho"] & (\iota_{\cB})^e_! \, P_{\cB/\cA}[n] \ar[r] & 0
	\end{tikzcd}
\end{equation}
In particular, notice that since the terms of \eqref{S_Pi_A_SES} are cofibrant, and hence in particular projective as graded modules, it remains exact after applying $F^e_!$ to it, so that the first row of \eqref{S_Pi_SES_A_to_BA} is indeed exact.

Taking the cone of the vertical maps of \eqref{S_Pi_SES_A_to_BA}, we see that there is a short exact sequence in $\Mod((\Pi_{\cB/\cA})^e)$:
\begin{equation*}
	0 \raq (\iota_{\cB})^e_! \, \cS(\cB,\cA) \xraq{\gamma} \cS(\Pi_{\cB/\cA}, \Pi_{\cA}) \xraq{\rho} (\iota_{\cB})^e_! \, \cone(\widetilde{\mu})[n] \raq 0
\end{equation*}

We will prove that $(\omega_2')^{\#}_L$ is a quasi-isomorphism by giving an analogue of Lemma \ref{triang_decomp_lemma} (see Lemma \ref{triang_decomp_lemma_rel}).
As a preparation for it, consider the canonical map
\begin{equation*}
	\begin{split}
		j \, : \, \Theta_{\cB/\cA} \, &:= \, \cone [ P_{\cA}[-1] \otimes_{\cAe} \cSA \xra{\mu \otimes \gamma_F} P_{\cB/\cA} \otimes_{\cBe} \cSB ] \\ 
		\, & \ra \, 
		\cone [ ((\iota_{\cA})^e_! \, P_{\cA}[-1]) \otimes_{(\Pi_{\cA})^e} ( (\iota_{\cA})^e_! \, \cSA) \xra{\mu \otimes \gamma_F} ((\iota_{\cB})^e_! \, P_{\cB/\cA}) \otimes_{(\Pi_{\cB/\cA})^e} ((\iota_{\cB})^e_! \, \cSB) ]
	\end{split}
\end{equation*}
and let $j(\theta_{\cB/\cA})$ be the image of the relative Casimir element $\theta_{\cB/\cA}$ under this map. 

As in the non-derived version of \eqref{rel_Hoch_induce_map_tensor} (see the proof of Proposition \ref{rel_Hoch_induce_map_of_triag_prop} and Remark \ref{non_cof_M_N_remark}), there are induced elements
\begin{equation*}
	\begin{split}
		j(\theta_{\cB/\cA})' \, &\in \, ((\iota_{\cB})^e_! \, P_{\cB/\cA}) \otimes_{(\Pi_{\cB/\cA})^e} ((\iota_{\cB})^e_! \, \cS(\cB,\cA)) \\
		j(\theta_{\cB/\cA})''' \, &\in \, ((\iota_{\cB})^e_! \, \cone(\widetilde{\mu})) \otimes_{(\Pi_{\cB/\cA})^e} ((\iota_{\cB})^e_! \, \cSB)
	\end{split}
\end{equation*}
which in turn induces maps in $\Mod((\Pi_{\cB/\cA})^e)$:
\begin{equation}  \label{j_theta_BA_prime_maps}
	\begin{split}
		(j(\theta_{\cB/\cA})')^{\#}_L \, &: \, ((\iota_{\cB})^e_! \, P_{\cB/\cA})^{\vee} \raq (\iota_{\cB})^e_! \, \cS(\cB,\cA) \\
		(j(\theta_{\cB/\cA})''')^{\#}_R \, &: \, ((\iota_{\cB})^e_! \, \cSB)^{\vee} \raq (\iota_{\cB})^e_! \, \cone(\widetilde{\mu})
	\end{split}
\end{equation}

\blm
The maps \eqref{j_theta_BA_prime_maps} are both quasi-isomorphisms.
\elm

\bpf
This proof is completely analogous to the proof of Lemma \ref{Casimir_element_lemma}(1),(4). First, we consider the maps in $\Mod(\cBe)$ induced by $\theta_{\cB/\cA}$, defined in the same way as \eqref{j_theta_BA_prime_maps} above:
\begin{equation}  \label{theta_BA_prime_maps}
	\begin{split}
		((\theta_{\cB/\cA})')^{\#}_L \, &: \, (P_{\cB/\cA})^{\vee} \raq  \cS(\cB,\cA) \\
		((\theta_{\cB/\cA})''')^{\#}_R \, &: \, (\cSB)^{\vee} \raq  \cone(\widetilde{\mu})
	\end{split}
\end{equation}

We claim that both of the maps \eqref{theta_BA_prime_maps} are quasi-isomorphisms. To show that $((\theta_{\cB/\cA})''')^{\#}_R$ is a quasi-isomorphism, notice that by the proof of Proposition \ref{rel_Hoch_induce_map_of_triag_prop} (see also Remark \ref{non_cof_M_N_remark}), there is a map of distinguished triangles in $\cD(\cBe)$:
\begin{equation*}
	\begin{tikzcd}
		\cSB^{\vee} \ar[d, "(\theta_{\cB/\cA}''')^{\#}_R"] \ar[r] & F^e_!(\cSA)^{\vee} \ar[d, "(\theta_{\cB/\cA}'')^{\#}_R"] \ar[r] & \cS(\cB,\cA)^{\vee}[1] \ar[d, "(\theta_{\cB/\cA}')^{\#}_R"] \\
		\cone(\widetilde{\mu}) \ar[r]  & F^e_! \, P_{\cA} \ar[r]  & P_{\cB/\cA}[1]
	\end{tikzcd}
\end{equation*}
The fact that $\theta_{\cB/\cA}$ is a relative Casimir element implies that $(\theta_{\cB/\cA}'')^{\#}_R$ and $(\theta_{\cB/\cA}')^{\#}_R$ are quasi-isomorphisms. Thus, $(\theta_{\cB/\cA}''')^{\#}_R$ is also a quasi-isomorphism.

To see that $((\theta_{\cB/\cA})')^{\#}_L$ is a quasi-isomorphism, use the diagram \eqref{PQ_theta_dual} for $Q = \cS(\cB,\cA)$, $P = P_{\cB/\cA}$ and $p = p_{\cB/\cA}$, and follow the argument in the paragraph that follows \eqref{PQ_theta_dual}.

From this, the fact that \eqref{j_theta_BA_prime_maps} are quasi-isomorphisms then follows the same proof as in Lemma \ref{Casimir_element_lemma}(4).
\epf

\blm \label{triang_decomp_lemma_rel}
Let $\zeta \in X^{(2)}(\cB,\cA)_n$ be an element such that $\omega_2 := B(J(s^{n-1} \theta_{\cB/\cA})) + I(\zeta) \in X^{(2)}(\Pi_{\cB/\cA},\Pi_{\cA})_n$ is $b$-closed, then the following diagram in $\Mod((\Pi_{\cB/\cA})^e)$ commutes:
\begin{equation}  \label{omega_triang_diag_rel}
	\begin{tikzcd}
		0 \ar[r] & ((\iota_{\cB})^e_! \, P_{\cB/\cA})^{\vee}  \ar[r,"\rho^{\vee}"] \ar[d, "(j(\theta_{\cB/\cA})')^{\#}_L"] & \cS(\Pi_{\cB/\cA})^{\vee}[n] \ar[r, "\gamma^{\vee}"] \ar[d, "(\omega_2')^{\#}_L"] & ((\iota_{\cB})^e_! \, \cSB)^{\vee}[n] \ar[r] \ar[d, "(j(\theta_{\cB/\cA})''')^{\#}_R"] & 0 \\
		0 \ar[r]  & (\iota_{\cB})^e_! \, \cS(\cB,\cA) \ar[r, "\gamma"] & \cS(\Pi_{\cB/\cA}, \Pi_{\cA}) \ar[r, "\rho"] & (\iota_{\cB})^e_! \, \cone(\widetilde{\mu})[n] \ar[r] & 0
	\end{tikzcd}
\end{equation}
where the maps $(j(\theta_{\cB/\cA})')^{\#}_L$ and $(j(\theta_{\cB/\cA})''')^{\#}_R$ are defined in \eqref{j_theta_BA_prime_maps}. 
\elm

\bpf
As in Lemma \ref{triang_decomp_lemma}, the pre-map $(I(\zeta)')^{\#}_L : \cS(\Pi_{\cB/\cA})^{\vee}[n] \ra \cS(\Pi_{\cB/\cA}, \Pi_{\cA})$ is the composition
\begin{equation*}
	\cS(\Pi_{\cB/\cA}, \Pi_{\cA}) \xraq{\gamma^{\vee}} ((\iota_{\cB})^e_! \, \cSB)^{\vee}[n] \xraq{(j(\zeta)')^{\#}_L}  (\iota_{\cB})^e_! \, \cS(\cB,\cA) \xraq{\gamma}  \cS(\Pi_{\cB/\cA}, \Pi_{\cA}) 
\end{equation*}
Hence, we may assume that $\zeta = 0$, and $\omega_2 := B(J(s^{n-1} \theta_{\cB/\cA}))$, which is not necessarily closed. Then we verify that the diagram \eqref{omega_triang_diag_rel} of pre-maps commutes for any element $\theta_{\cB/\cA} \in \Theta(\cB,\cA)_0$, not necessarily a relative Casimir element, not even necessarily closed.

For this, we consider $\theta_{\cB/\cA} \in \Theta(\cB,\cA)_0$ to be one of the following $4$ types of elements:
\begin{enumerate}
	\item[(A1)] $s(\xi_{\cA} \otimes sDf)$ for $f \in \cA(x,y)_{-m-2}$ and  $\xi_{\cA} \in P_{\cA}(y,x)_m$
	\item[(A2)] $s(\xi_{\cA} \otimes E_x^{\cA})$ for $\xi_{\cA} \in P_{\cA}(x,x)_{-1}$
	\item[(B1)] $\xi_{\cB} \otimes sDg$ for $g \in \cB(x,y)_{-m-1}$ and  $\xi_{\cB} \in P_{\cB/\cA}(y,x)_m$
	\item[(B2)] $\xi_{\cB} \otimes E_x^{\cB}$ for $\xi_{\cB} \in P_{\cB/\cA}(x,x)_{0}$
\end{enumerate}
Then by Lemma \ref{B_on_X1}, the element $\omega_2 = B(J(s^{n-1} \theta_{\cB/\cA})) \in \cone[X^{(2)}(\Pi_{\cA}) \ra X^{(2)}(\Pi_{\cB/\cA})]$ is given by 
\begin{enumerate}
	\item[(A1)] $s(sD(s^{n-2}\xi_{\cA}) \otimes sDf + (-1)^{(m+n-1)(m+1)} sDf \otimes sD(s^{n-2}\xi_{\cA}))  \, \in X^{(2)}(\Pi_{\cA})[1]$
	\item[(A2)] $s(sD(s^{n-2}\xi_{\cA}) \otimes E_x^{\cA} + E_x^{\cA} \otimes sD(s^{n-2}\xi_{\cA})) \, \in X^{(2)}(\Pi_{\cA})[1]$
	\item[(B1)] $sD(s^{n-1}\xi_{\cB}) \otimes sDg + (-1)^{m(m+n)} sDg \otimes sD(s^{n-1}\xi_{\cB}) \, \in X^{(2)}(\Pi_{\cB/\cA})$
	\item[(B2)] $sD(s^{n-1}\xi_{\cB}) \otimes E_x^{\cB} + E_x^{\cB} \otimes sD(s^{n-1}\xi_{\cB}) \, \in X^{(2)}(\Pi_{\cB/\cA})$
\end{enumerate}
Thus, the element $\omega_2' \in \cS(\Pi_{\cB/\cA}) \otimes_{\cBe} \cone[ F^e_!(\cS(\Pi_{\cA})) \ra \cS(\Pi_{\cB/\cA}) ]$ is given by
\begin{enumerate}
	\item[(A1)] $sD(s^{n-2}\mu(\xi_{\cA})) \otimes s(sDf) + (-1)^{(m+n-1)(m+1)} sD(F(f)) \otimes s(sD(s^{n-2}\xi_{\cA}))$
	\item[(A2)] $sD(s^{n-2}\mu(\xi_{\cA})) \otimes sE_x^{\cA} + E_x^{\cB} \otimes s(sD(s^{n-2}\xi_{\cA})))$
	\item[(B1)] $sD(s^{n-1}\xi_{\cB}) \otimes sDg + (-1)^{m(m+n)} sDg \otimes sD(s^{n-1}\xi_{\cB})$
	\item[(B2)] $sD(s^{n-1}\xi_{\cB}) \otimes E_x^{\cB} + E_x^{\cB} \otimes sD(s^{n-1}\xi_{\cB})$
\end{enumerate}
where the ones for type (A1), (A2) are in $\cS(\Pi_{\cB/\cA}) \otimes_{\cBe} ( F^e_!(\cS(\Pi_{\cA}))[1])$, while the ones for (B1), (B2) are in $\cS(\Pi_{\cB/\cA}) \otimes_{\cBe} \cS(\Pi_{\cB/\cA})$. 
For each of these types, one can directly verify that \eqref{omega_triang_diag_rel} commutes.
\epf

This completes the proof that $\widetilde{\omega}$ is a relative $n$-Calabi-Yau structure. The proof of the statement about exactness of  $\widetilde{\omega}$ in the absolute case can be carried over almost verbatim to the present relative case.
\epf

\end{document}